\documentclass[11pt, final]{amsart}

\title{On the $K$-theory of coordinate axes in affine space}
\author{Martin Speirs}

\usepackage{import}
\usepackage{ThesisPreamble}
\usepackage{verbatim}
\usepackage{soul}
\usepackage{booktabs}
\usepackage[square,sort,comma,numbers]{natbib}

\newcommand{\HH}{\operatorname{HH}}
\newcommand{\HP}{\operatorname{HP}}
\newcommand{\THH}{\operatorname{THH}}

\newcommand{\TP}{\operatorname{TP}}

\newcommand{\TC}{\operatorname{TC}}
\newcommand{\HC}{\operatorname{HC}}
\newcommand{\cyc}{\operatorname{cyc}}
\newcommand{\can}{\operatorname{can}}
\newcommand{\Spec}{\operatorname{Spec}}

\newcommand{\Bcy}{\operatorname{B^{cy}}}
\newcommand{\Fil}{\operatorname{Fil}}
\newcommand{\sd}{\operatorname{sd}}
\newcommand{\Htilde}{\tilde{H}}

\begin{document}
\maketitle

\section{Introduction}
\noindent
Let $k$ be a perfect field of characteristic $p > 0$ and let $A_d$ denote the $k$-algebra $k[x_1, \dots , x_d]/(x_ix_j)_{i \ne j}$, which is the coordinate ring of the coordinate $(x_1, \dots , x_d)$-axes in affine $d$-space $\Aff^{d} = \Spec(k[x_1, \dots, x_d])$ over $k$. This is an affine curve with a singularity at the origin, the singularity being defined by the ideal $I_d = (x_1, \dots , x_d)$. Our main result is the computation of the relative algebraic $K$-theory, $K(A_d, I_d)$ of the pair $(A_d,I_d)$. The relative $K$-theory is defined to be the mapping fiber of the map $K(A_d) \to K(k)$ induced by the projection onto $k = A_d/I_d$.  
% When $d > 2$ this ring is no longer a complete intersection. 

To state the result we introduce some notation. We consider words in $d$ letters $x_1, \dots , x_d$, i.e.\ a finite string $\omega = w_1w_2 \dots w_m$ where each $w_i$ is one of the letters $x_1, \dots , x_d$. A word $\omega = w_1w_2 \dots w_m$ has \emph{no cyclic repetitions} if $w_i \ne w_{i+1}$ for all $i = 1, \dots , m-1$ and if $w_m \ne w_1$. 
For $d \geq 1$ and $s \geq 1$ let $\cyc_d(s)$ denote the number of cyclic words in $d$ letters, of length $s$, period $s$, having no cyclic repetitions. In \cref{CyclicAppendix} we give a formula for $\cyc_d(s)$.  Suppose $m' \geq 1$ and let $t_{ev} = t_{ev}(p,r,m')$ be the unique positive integer such that $p^{t_{ev}-1}m' \leq 2r < p^{t_{ev}}m'$, or zero if such $t_{ev}$ does not exist.  
Let $t_{od} = t_{od}(p,r,m')$ be the unique positive integer such that $p^{t_{od}-1}m' \leq 2r + 1 < p^{t_{od}}m'$, or zero if such $t_{od}$ does not exist. Let $J_p$ denote the set of positive integers $m' \geq 1$ which are coprime to $p$.

\begin{theorem}\label{KtheoryCoordinateAxes}
Let $k$ be a perfect field of characteristic $p >0$. Let $A_d$ be the ring $k[x_1, \dots , x_d]/(x_ix_j)_{i \ne j}$ of coordinate axes, and let $I_d = (x_1, \dots , x_d)$. 
Then if $p > 2$
\[
K_q(A_d,I_d)  \cong   
   \left\{
        \begin{array}{ll}
            \prod_{\substack{m' \in J_p \\ \text{even}}} \prod_{\substack{s' \mid m' \\ \text{even}}} \prod_{u \leq t_{ev}} W_{t_{ev}-u}(k)^{\oplus \cyc_d(p^us')} & \quad q = 2r \\
            \prod_{\substack{m' \in J_p \\ \text{odd}}} \prod_{s' \mid m'} \prod_{u \leq t_{od}} W_{t_{od}-u}(k)^{\oplus \cyc_d(p^us')}  & \quad  q = 2r + 1
        \end{array}
    \right.
\]
 If $p = 2$ then
\[
K_q(A_d,I_d)  \cong   
    \left\{
        \begin{array}{ll}
            \prod_{\substack{m' \geq 1 \\ \text{odd}}} \prod_{s' \mid m'} \prod_{1 \leq u \leq t_{ev}} W_{t_{ev}-u}(k)^{\oplus \cyc_d(2^us')} & \quad q = 2r \\
           \prod_{\substack{m' \geq 1 \\ \text{odd}}} \prod_{s' \mid m'} \prod_{0 \leq v \leq t_{ev}} k^{\oplus \cyc_d(s')}  & \quad  q = 2r + 1 
        \end{array}
    \right.
\]
In both cases $t_{ev} = t_{ev}(p,r,m')$ and $t_{od} = t_{od}(p,r,m')$ are as defined above.
\end{theorem}

\noindent Note that the products appearing in the statement are all finite since for $m'$ large enough $t_{ev} = t_{od} = 0$.

The result extends the calculation by Dennis and Krusemeyer when $q=2$, \cite[Theorem 4.9]{DennisKrusemeyer}. Hesselholt carried out the computation in the case $d=2$ in \cite{Hesselholt2007}. Our strategy follows the one in \cite{Hesselholt2007}, except that we use the framework for $\TC$ set up by Nikolaus and Scholze \cite{NikolausScholze} to which we refer for background on $\TC$ and cyclotomic spectra. 

The analogous computation when working with a base field of characteristic zero was carried out by Geller, Reid and Weibel in \cite[Theorem 7.1.]{GellerReidWeibel}. In \cref{CharZeroSection} we provide an extension of their calculations to the case of ind-smooth $\QQ$-algebras.
 
Our proof of \cref{KtheoryCoordinateAxes} uses the cyclotomic trace map from $K$-theory to topological cyclic homology, $\TC$. See also \cite{HesselholtMadsen97} for background on cyclotomic spectra and for similar calculations. We remark that Hesselholt and Nikolaus \cite{HesselholtNikolausCusp} have completed a calculation for $K$-theory of cuspical curves using similar methods as this paper.

\subsection{Overview}
In \cref{KABI} we reduce the $K$-theory computation to a $\TC$ computation. \cref{THHsection} carries out the requisite $\THH$ computation. Then \cref{TPsection} and \cref{TCsection} assemble this to complete the proof. \cref{TPsection} contains a new method for computing $\TP$, which makes crucial use of the Nikolaus-Scholze framework, see also \cite{SpeirsTruncated} for a similar calculation. In \cref{CharZeroSection} we consider the characteristic zero situation and extend the computation of \cite{GellerReidWeibel}. Finally in \cref{CyclicAppendix} we derived the necessary counting formula for cyclic words.

\subsection{Acknowledgements}
I am grateful to my advisor Lars Hesselholt for his guidance and support during the production of this paper. Special thanks are due to Fabien Pazuki for encouraging and useful conversations. I thank Benjamin B\"ohme, Ryo Horiuchi, Joshua Hunt, Mikala Jansen, Manuel Krannich, and Malte Leip for several useful conversations. It is a pleasure to thank Malte Leip for carefully reading a draft version of this paper and providing several corrections and suggestions. I thank the anonymous referee for a giving the paper a detailed reading, spotting several mistakes, and providing many helpful suggestions. This paper is based upon work supported by the DNRF Niels Bohr Professorship of Lars Hesselholt
as well as the National Science Foundation under Grant No.\ DMS-1440140 while the author was in residence at the Mathematical Sciences Research Institute in Berkeley, California during the Spring 2019 semester. Additionally part of the work was supported by the Independent Research Fund Denmark, Grant No.\ 9056-00018B with mentorship from Martin Olsson.

%%%%%%%%%%%%%%%%%%%%%%%%%%%%%%%%%%%%%%%%
%%%%%%%%%%%%%%%%%%%%%%%%%%%%%%%%%%%%%%%%
%%%%%%%%%%%%%%%%%%%%%%%%%%%%%%%%%%%%%%%%

\section{Bi-relative $K$-theory and topological cyclic homology}\label{KABI}
\noindent
The normalization of $A_d$ is just $d$ disjoint lines, whose coordinate ring is 
\[
B_d = k[x_1] \times \cdots \times k[x_d]. 
\]
Gluing these lines together at $x_1 = x_2 =  \dots = x_d = 0$ one obtains $\Spec(A_d)$. Algebraically this is the statement that the following square is a pullback of rings
\[
\begin{tikzcd}[row sep=3em, column sep=3em]
A_d \arrow[r] \arrow[d] & k \arrow[d] \\
B_d \arrow[r] & k^{\times d}
\end{tikzcd}
\]
Here the horizontal maps take the variables $x_i$ to zero. The left-vertical map is the normalization map. It maps $x_i$ to $(0, \dots, x_i \dots, 0)$ where $x_i$ is in the $i$'th position. The right-vertical map is the diagonal. If $K$-theory preserved pullbacks then the diagram would give a computation of $K(A_d)$ in terms of $K(B_d)$ and $K(k)$. Using the fundamental theorem of $K$-theory (since $k$ is regular) one would get a formula for $K(A_d)$ purely in terms of $K(k)$. But $K$-theory does not preserve pullbacks. However, there is still something to be done. We can form the bi-relative $K$-theory, $K(A_d,B_d, I_d)$ as the iterated mapping fiber of the diagram
\[\begin{tikzcd}[row sep=3em, column sep=3em]
K(A_d) \arrow[r] \arrow[d] & K(k) \arrow[d] \\
K(B_d) \arrow[r] & K(k^{\times d})
\end{tikzcd}\]
Again, if $K$-theory preserved pullbacks then $K(A_d,B_d,I_d)$ would be trivial, but (as we shall see) it is not. Since $k$ is regular the fundamental theorem of algebraic $K$-theory \cite[Section 6]{Quillen73}, and the fact that $K$-theory \emph{does} preserve products, shows that the canonical map
\[
K(A_d,B_d,I_d) \longto K(A_d,I_d)
\]
is an equivalence. Here $K(A_d,I_d)$ is the relative $K$-theory spectrum, i.e.\ the mapping fiber of the map $K(A_d) \to K(A_d/I_d)$ induced by the quotient. 
Any splitting of the quotient map $A_d \to A_d/I_d$ provides a splitting of $K$-groups
\[
K_q(A_d) \cong K_q(k) \oplus K_q(A_d,I_d).
\]
So if we can compute $K_q(A_d,I_d) = K_q(A_d,B_d,I_d)$ we have a computation of $K_q(A_d)$. This is what we will do.

Geisser and Hesselholt \cite{GeisserHesselholtBiRelative} have shown that the cyclotomic trace induces an isomorphism
\[
K_q(A,B,I, \ZZ/p^\nu) \overset{\sim}\longto \TC_q(A,B,I;p, \ZZ/p^\nu)
\]
between the bi-relative $K$-theory and the bi-relative topological cyclic homology, when both are considered with $\ZZ/p^v$ coefficients. 
The corresponding statement with rational coefficients (using Connes-Tsygan negative cyclic homology and the Chern character) was proven by Corti\~nas in \cite{CortinasKABI}.

In fact it suffices to compute $\TC(A_d,B_d,I_d, p; \ZZ_p)$ since the trace map 
\[
K(A,B,I) \to \TC(A,B,I, p; \ZZ_p)
\]
 is an equivalence whenever $p$ is nilpotent in $A$, as shown in \cite[Theorem C]{GeisserHesselholt2011}. Furthermore, since $\THH(A_d,B_d,I_d;p)$ is an $Hk$-module, it is in particular $p$-complete, and so $\TC(A,B,I, p ; \ZZ_p) \simeq \TC(A,B,I,p)$ (see \cite[Section II.4]{NikolausScholze}).
Thus, to prove \cref{KtheoryCoordinateAxes} it suffices to prove the following result.

\begin{theorem}\label{TCCoordinateAxes}
Let $k$ be a perfect field of characteristic $p >0$. Let  $A_d$ be the ring $k[x_1, \dots , x_d]/(x_ix_j)_{i \ne j}$ of coordinate axes, $B_d = k[x_1] \times \dots \times k[x_d]$ and let $I_d = (x_1, \dots , x_d)$.
Then if $p > 2$
\begin{eqnarray*}
\TC_q(A_d, B_d, I_d)  \cong   
    \left\{
        \begin{array}{ll}
            \prod_{\substack{m' \in J_p \\ \text{even}}} \prod_{\substack{s' \mid m' \\ \text{even}}} \prod_{u \leq t_{ev}} W_{t_{ev}-u}(k)^{\oplus \cyc_d(p^us')} & \quad q = 2r \\
            \prod_{\substack{m' \in J_p \\ \text{odd}}} \prod_{s' \mid m'} \prod_{u \leq t_{od}} W_{t_{od}-u}(k)^{\oplus \cyc_d(p^us')}  & \quad  q = 2r {+} 1
        \end{array}
    \right.
\end{eqnarray*}
If $p = 2$ then
\begin{eqnarray*}
\TC_q(A_d, B_d, I_d)  \cong   
    \left\{
        \begin{array}{ll}
            \prod_{\substack{m' \geq 1 \\ \text{odd}}} \prod_{s' \mid m'} \prod_{1 \leq u \leq t_{ev}} W_{t_{ev}-u}(k)^{\oplus \cyc_d(2^us')} & \quad q = 2r \\
           \prod_{\substack{m' \geq 1 \\ \text{odd}}} \prod_{s' \mid m'} \prod_{0 \leq v \leq t_{ev}} k^{\oplus \cyc_d(s')}  & \quad  q = 2r {+} 1  .
        \end{array}
    \right.
\end{eqnarray*}
\end{theorem}

%%%%%%%%%%%%%%%%%%%%%%%%%%%%%%%%%%%%%%%%
%%%%%%%%%%%%%%%%%%%%%%%%%%%%%%%%%%%%%%%%
%%%%%%%%%%%%%%%%%%%%%%%%%%%%%%%%%%%%%%%%

\section{$\THH$ and the cyclic bar construction}\label{THHsection}

\noindent The spectrum $\TC(A_d, B_d, I_d)$ is defined using topological Hochschild homology, so that is where we start. In this section we drop the subscript $d$, so that $A = A_d$, $B = B_d$ and $I = I_d$.  The bi-relative topological Hochschild homology is the spectrum $\THH(A,B,I)$ defined as the iterated mapping fiber of the following diagram.
\begin{equation}\label{THHdiagram}
\begin{tikzcd}[row sep=3em, column sep=3em]
\THH(A) \arrow[r] \arrow[d] & \THH(A/I) \arrow[d] \\
\THH(B) \arrow[r] & \THH(B/I)
\end{tikzcd}
\end{equation}

The ring $A$ is a pointed monoid ring and, by \cref{SplittingLemma} below, we may compute $\THH(A)$  in terms of the cyclic bar construction of the defining pointed monoid. 

\subsection{Unstable cyclotomic structure on the cyclic bar construction} \label{UnstableFrobenius}% \todo{formulate a lemma to the effect that the unstable Frobenius on cyclic bar of monoid is an equivalence. This is used later in the lemma on connectivity.}
\noindent Let $\Pi$ be a pointed monoid, that is a monoid object in the symmetric monoidal category of based spaces and smash product. The cyclic bar construction of $\Pi$ is the cyclic space $\Bcy(\Pi)[-]$ with
\[
\Bcy(\Pi)[k] = \Pi^{\wedge (k+1)}
\]
and with the usual Hochschild-type structure maps.
\begin{eqnarray*}
d_i(\pi_0 \wedge \dots \wedge \pi_m) &=& \left\{
        \begin{array}{ll}
         \pi_0 \wedge \dots \wedge \pi_i \pi_{i+1} \wedge \dots \wedge \pi_{m}    & \quad 0 \leq i < m \\
          \pi_{m}\pi_0 \wedge \pi_1 \wedge \dots  \wedge \pi_{m-1} & \quad  i = m
        \end{array}
    \right.
\\
s_i(\pi_0 \wedge \dots \wedge \pi_m) &=&  \pi_0 \wedge \dots \wedge \pi_i \wedge 1 \wedge \pi_{i+1} \wedge \dots \wedge \pi_m
\\
t_m(\pi_0 \wedge \dots \wedge \pi_m) &=& \pi_m \wedge \pi_0 \wedge \dots \wedge \pi_{m-1}
\end{eqnarray*}

 We write $\Bcy(\Pi)$ for the geometric realization of $\Bcy(\Pi)[-]$. It is a space with $\TT$-action where $\TT$ is the circle group (for a proof see for example \cite[Theorem 7.1.4.]{LodayBook}). Furthermore it is an unstable cyclotomic space, %in the sense of \cite{NikolausCyclicK},
  i.e.\ there is a map
\[
\psi_p : \Bcy(\Pi) \to \Bcy(\Pi)^{C_p}
\] 
which is equivariant when the codomain is given the natural $\TT/C_p$-action. This map goes back to \cite[Section 2]{BHM}. We briefly sketch the construction. The $C_p$-action on $\Bcy(\Pi)$ is not simplicial, but we can make it so by using the edge-wise subdivision functor $\sd_p : \Delta \to \Delta$ which is given by the $p$-fold concatenation, $\sd_p [m-1]  = [m-1] \amalg \dots \amalg [m-1]$ and $\sd_p(\theta) = \theta \amalg \dots \amalg \theta$ for morphisms $\theta : [m-1] \to [n-1]$. Given a simplicial set $X[-]$ we let $\sd_p X[-] = X[-] \circ \sd^{op}_p$. For the topological simplex $\Delta^{m-1} \subset \RR^{m}$ let $d_p : \Delta^{m-1} \to \Delta^{pm-1}$ be the diagonal embedding
\[
d_p(z) = \frac{1}{p} z \oplus \dots \oplus \frac{1}{p} z .
\]
This induces a (non-simplicial) map on geometric realization 
\[
D_p : |(\sd_p X)[-]| \to | X[-]|
\]
 by $id \times d_p : X[pm-1] \times \Delta^{m-1} \to X[pm-1] \times \Delta^{pm-1}$. Then \cite[Lemma 1.1.]{BHM} shows that $D_p$ is a homeomorphism -- one need only check on the representables $\Delta^k[-]$ where it follows by explicit calculations. Furthermore, when $X[-]$ is a cyclic space, then $D_p$ is a $\TT$-equivariant homoemorphism \cite[Lemma 1.11.]{BHM}.  In the case $X[-] = \Bcy(\Pi)[-]$ one has available the simplicial diagonal 
 \[
 \tilde{\Delta_p} : \Bcy(\Pi)[- ] \to \sd_p\Bcy(\Pi)[-]
 \]
given by
\[
\overrightarrow{\pi} = \pi_0 \wedge \dots  \wedge \pi_m \quad \longmapsto \quad \overrightarrow{\pi} \wedge \dots \wedge \overrightarrow{\pi} \qquad \text{($p$ copies)}
\]
It clearly lands in the $C_p$-fixed points of $\sd_p\Bcy(\Pi)[-]$ and induces an isomorphism $\tilde{\Delta_p} : \Bcy(\Pi)[- ] \to \sd_p\Bcy(\Pi)[-]^{C_p}$. We define 
\[
\psi_p : \Bcy(\Pi) \to \Bcy(\Pi)^{C_p}
\]
 to be the composite
\begin{eqnarray*}
\Bcy(\Pi) = | \Bcy(\Pi)[-]| &\underset{\cong}{\overset{\tilde{\Delta_p}}\longto} & | \sd_p \Bcy(\Pi)[-]^{C_p}|
\\ &\underset{\cong}{\longto} & |\sd_p \Bcy(\Pi)[-]|^{C_p} 
\\ &\underset{\cong}{\overset{D_p}\longto} & |\Bcy(\Pi)[-]|^{C_p} = \Bcy(\Pi)^{C_p}
\end{eqnarray*}
where the middle map is the canonical equivalence witnessing the fact that geometric realization commutes with finite limits. As all three maps in the definition are $\TT$-equivariant, so is $\psi_p$. 

Passing to suspension spectra now gives us a cyclotomic structure on $\THH(\SSS(\Pi))$ where $\SSS(\Pi) = \Sigma^\infty\Pi$. Indeed, there is always a $\TT/C_p$-equivariant map $\Bcy(\Pi)^{C_p} \to \Bcy(\Pi)^{hC_p}$. Composing with this map gives a $\TT \simeq \TT/C_p$-equivariant map $\psi^h : \Bcy(\Pi) \to \Bcy(\Pi)^{hC_p}$. Thus, we obtain a $\TT \simeq \TT/C_p$-equivariant map
\begin{eqnarray*}
\THH(\SSS(\Pi)) = \Sigma^\infty \Bcy(\Pi) &\overset{\psi^h}\longto & \Sigma^\infty \Bcy(\Pi)^{hC_p} 
\\ &\longto & (\Sigma^\infty \Bcy(\Pi))^{hC_p} 
\\ &\overset{can}\longto & (\Sigma^\infty \Bcy(\Pi))^{tC_p} = \THH(\SSS(\Pi))^{tC_p}
\end{eqnarray*}
where the middle map is the canonical map from the suspension spectrum of a homotopy limit, to the homotopy limit of the suspension spectrum and  the last map is the canonical map from the homotopy fixed points to the Tate construction.

\begin{comment}
The following is a pointed version of Lemma IV. 3.1. in \cite{NikolausScholze}.
\begin{lemma} \todo{replace this lemma with the one using actual fixed points.}
Let $\Pi$ be a pointed monoid in spaces.
\begin{enumerate}
\item The cyclic bar construction $\Bcy(\Pi)$ admits a canonical $\TT$-equivariant map
\[
\psi_p : \Bcy(\Pi) \to \Bcy(\Pi)^{hC_p}
\]
for every prime $p$.
\item There is a commutative diagram
\begin{equation*}
\begin{tikzcd}[row sep=2em, column sep=2em]
\Pi \arrow[r, "i"] \arrow[d, "\Delta"] & \Bcy(\Pi) \arrow[d, "\psi_p"] \\
(\Pi \wedge \dots \wedge \Pi)^{hC_p} \arrow[r] & \Bcy(\Pi)^{hC_p}
\end{tikzcd}
\end{equation*}
\item Upon taking suspension spectrum the map $\psi_p$ refines the cyclotomic structure of $\THH(\SSS(\Pi)) = \Sigma^\infty \Bcy(\Pi)$ in the sense that there is a commutative diagram
\begin{equation*}
\begin{tikzcd}[row sep=3em, column sep=3em]
\Sigma^\infty \Bcy(\Pi) \arrow[r, "\psi_p"] \arrow[d, "\phi_p"] & \Sigma^\infty (\Bcy(\Pi)^{hC_p}) \arrow[d] \\
(\Sigma^\infty \Bcy(\Pi))^{tC_p} & (\Sigma^\infty \Bcy(\Pi))^{hC_p}  \arrow[l , "can"]
\end{tikzcd}
\end{equation*}
where $\phi_p$ is the cyclotomic Frobenius map on $\THH(\SSS(\Pi))$.
\end{enumerate}
\end{lemma}
\end{comment}

\subsection{$\THH$ of monoid algebras}

Both of the algebras $A$ and $B$ (and their quotients by $I$) are pointed monoid algebras, enabling us to use the following splitting result.

\begin{lemma}{\rm \cite[Theorem 7.1.]{HesselholtMadsen97}} \label{SplittingLemma}
Let $k$ be a ring, $\Pi$ a pointed monoid, and $k(\Pi)$ the pointed monoid algebra. Let $\Bcy(\Pi)$ be the cyclic bar-construction on $\Pi$. Then there is a $\TT$-equivariant equivalence
\[
\THH(k) \otimes \Bcy(\Pi) \overset{\sim}\to \THH(k(\Pi)).
\]
Under this equivalence, the Frobenius on $\THH(k(\Pi))$ is identified with the tensor product of the Frobenius on $\THH(k)$ with the one on $\THH(\SSS(\Pi))$, and as such is an equivalence of cyclotomic spectra.
\end{lemma}
\begin{proof}
Since $\THH$ is symmetric monoidal  \cite[Section IV.2.]{NikolausScholze} we obtain
\[
\THH(k(\Pi)) = \THH(k \otimes \SSS(\Pi)) \simeq \THH(k) \otimes \THH(\SSS(\Pi)).
\]
Since  $\THH(\SSS(\Pi)) = \SSS(\Bcy(\Pi))$ we obtain 
\[
\THH(k(\Pi)) \simeq \THH(k) \otimes \SSS(\Bcy(\Pi))
\] as claimed.
\end{proof}

Let $\Pi^d = \{0, 1, x_1, x_1^2 , \dots , x_2 , x_2^2 , \dots , x_d , x_d^2 , \dots \}$ be the multiplicative pointed monoid with basepoint $0$ and multiplication determined by $x_ix_j = 0$ when $i \ne j$. Then $A_d \cong k(\Pi^d)$ is the pointed monoid algebra on $\Pi^d$. Let $\Pi^1 = \{0,1, t , t^2, \dots\}$ and $\Pi^0 = \{0, 1 \}$, so $B = k( \Pi^1 ) \times \cdots \times k(\Pi^1)$ (with $d$ products) and $k = k(\Pi^0)$. The diagram of cyclotomic spectra (\ref{THHdiagram}) is induced by the diagram of pointed monoids
\begin{equation}\label{Monoiddiagram}
\begin{tikzcd}[row sep=4em, column sep=4em]
\Pi^d \arrow[r, "\epsilon"] \arrow[d, "{(\rho_1,\ \dots ,\ \rho_d)}"] & \Pi^0 \arrow[d, "\Delta"] \\
\Pi^1 \times \cdots \times \Pi^1 \arrow[r, "{\epsilon^{\times d} }"] & \Pi^0 \times \cdots \times \Pi^0
\end{tikzcd}
\end{equation}
The map $\rho_j$ takes $x_i$ and to $0$ when $i \ne j$ and takes $x_j$ to $t$. The map $\epsilon$ takes the variables $x_1, \dots, x_d$ and $t$ to $0$. The map $\Delta$ is the diagonal.

The cyclic bar-construction sometimes decomposes, as a pointed $\TT$-space, into a wedge of spaces simple enough that one can understand their $\TT$-homotopy type. To do this for $\Bcy(\Pi^d)$ we need the notion of cyclic words.

We consider the alphabet $S = \{x_1, \dots, x_d\}$ and words 
\[
\omega : \{1,2, \dots , m\} \to S.
\]
Here $m$ is the \emph{length} of $\omega$. The cyclic group $C_m$ acts on the set of words of length $m$. An orbit for this action is called a \emph{cyclic word}. Such a cyclic word $\overline{\omega}$ has a \emph{period} namely the cardinality of the orbit. Words may be concatenated to give new, longer words, though concatenation is not well-defined for cyclic words. The empty word $\emptyset \to S$ is the unit for concatenation. It has length $0$ and period $1$. 

We can associate words to non-zero elements of $\Bcy(\Pi^d)[m]$ as follows. If $\pi \in \Pi^d = \Bcy(\Pi^d)[0]$ is non-zero then it is  of the form $\pi = x_j^l$ for some $1 \leq j \leq d$ and $l \geq 0$. Let $\omega(\pi)$ be the unique word of length $l$ all of whose letters are $x_j$. For example $\omega(x_1^2) = x_1x_1$ and $\omega(1) = \emptyset$. Now for a non-zero element $\pi_0 \wedge \dots \wedge \pi_m \in \Bcy(\Pi^d)[m]$ we let 
\[
\omega(\pi_0 \wedge \dots \wedge \pi_m) = \omega(\pi_0) \star \dots \star \omega(\pi_m)
\]
be the concatenation of each of the words $\omega(\pi_j)$. Note that $\omega(1)$ is the empty word, which is the unit for concatenation. For a cyclic word $\overline{\omega}$ we define
\[
\Bcy(\Pi^d, \overline{\omega}) [m] \subs \Bcy(\Pi^d) [m]
\]
to be the subset consisting of the basepoint and all elements $\pi_0 \wedge \dots \wedge \pi_m$ such that 
\[
\omega(\pi_0 \wedge \dots \wedge \pi_m) \in \overline{\omega}.
\]
The cyclic structure maps preserve this property and so, as $m \geq 0$ varies, this defines a cyclic subset 
\[
\Bcy(\Pi^d, \overline{\omega})[-] \subs \Bcy(\Pi^d)[-].
\]
Denote by $\Bcy(\Pi^d, \overline{\omega})$ the geometric realization of this subset. We will also often abbreviate $B(\overline{\omega}) = \Bcy(\Pi^d, \overline{\omega})$. As $\overline{\omega}$ ranges over all cyclic words every non-zero $m$-simplex $\pi_0 \wedge \dots \wedge \pi_m$ appears in exactly one such cyclic subset $\Bcy(\Pi^d, \overline{\omega})$. Thus we get a decomposition
\[
\Bcy(\Pi^d) = \bigvee \Bcy(\Pi^d, \overline{\omega})
\]
indexed on the set of all cyclic words with letters in $S = \{x_1 ,\dots, x_d\}$. 

\begin{lemma}\label{IndexLemma}
There is a canonical $\TT$-equivariant equivalence 
\[
\bigoplus \THH(k) \otimes \Bcy(\Pi^d, \overline{\omega}) \overset{\sim}{\to} \THH(A,B,I)
\]
where the sum on the left-hand side is indexed over all cyclic words whose period is $\geq 2$.
\end{lemma}
\begin{proof}
We first note that $\THH$ preserves product of rings, in the sense that the canonical map $\THH(R \times R') \to \THH(R) \times \THH(R')$ is a $\TT$-equivariant equivalence \cite[Proposition 4.2.4.4.]{DGM}. Using these equivalences we may consider the diagram induced from (\ref{THHdiagram}) using \cref{SplittingLemma},
\begin{equation}\label{THHdiagram2}
\begin{tikzcd}[row sep=4em, column sep=3em]
\THH(k) \otimes \Bcy(\Pi^d) \arrow[r, "\epsilon"] \arrow[d, "{(\rho_1 ,\ \dots ,\ \rho_d)}"] 
& \THH(k) \otimes \Bcy(\Pi^0) \arrow[d, "\Delta"] \\
(\THH(k) \otimes \Bcy(\Pi^1))^{\times d}    \arrow[r, "\epsilon^{\times d}"] 
& (\THH(k) \otimes \Bcy(\Pi^0))^{\times d}
\end{tikzcd}
\end{equation}
$\THH(A,B,I)$ is the iterated mapping fiber of this diagram. The mapping fiber of the left-hand vertical map consists of two parts; the part of the wedge sum indexed on cyclic words containing at least two different letters from $S$ and a part on which $\epsilon$ is an equivalence.
The map $\epsilon$ is trivial on this part indexed on cyclic words containing at least two different letters, finishing the claim.
\end{proof}

%%%%%%%%%%%%%%%%%%%%%%%%%%%%%%%%%%%%%%%
%%%%%%%%%%%%%%%%%%%%%%%%%%%%%%%%%%%%%%%

\subsection{Homotopy type of $\Bcy(\Pi^d; \overline{\omega})$.}

In this section we determine the homotopy type of the subspaces $\Bcy(\Pi^d, \overline{\omega}) \subs \Bcy(\Pi^d)$.

\begin{definition}\label{NoReps}
We say that a word $\omega = w_1w_2 \dots w_m$ has \emph{no cyclic repetitions} if $w_i \ne w_{i+1}$ for all $i = 1, \dots , m-1$ and if $w_m \ne w_1$. If on the other hand this is not satisfied, we say $\omega$ (or $\overline{\omega}$) \emph{has cyclic repetitions}. 
\end{definition}

\begin{lemma}\label{ClassificationLemma}
Let $\overline{\omega}$ be a cyclic word of period $s \geq 2$, with letters in the alphabet $S = \{x_1, \dots , x_d\}$. The homotopy type of the pointed $\TT$-space $\Bcy(\Pi^d, \overline{\omega})$ is given as follows.
\begin{enumerate}
\item If $\overline{\omega}$ has length $m = si$ and has no cyclic repetitions then a choice of representative word $\omega \in \overline{\omega}$ determines a $\TT$-equivariant homeomorphism
\[
S^{\RR[C_m]-1} \wedge_{C_i} \TT_+ \overset{\sim}\to \Bcy(\Pi^d, \overline{\omega}),
\]
where $\RR[C_m]-1$ is the reduced regular representation of $C_m$.
\item If $\overline{\omega}$ has cyclic repetitions, then $\Bcy(\Pi^d,\overline{\omega})$ is $\TT$-equivariantly contractible.
\end{enumerate}
\end{lemma}
\begin{proof}
This proof follows closely that of \cite[Lemma 1.6]{Hesselholt2007}. 
Let $\overline{\omega}$ be a cyclic word of length $m = si$ with period $s \geq 2$. Choose $\omega \in \overline{\omega}$ and let $(\pi_0, \dots , \pi_{m-1})$ be the unique $m$-tuple of non-zero elements in $\Pi^d$ such that $\omega(\pi_0 \wedge \dots  \wedge \pi_{m-1}) = \omega$. The pointed cyclic set $\Bcy(\Pi^d, \overline{\omega})[-]$ is generated by the $(m-1)$-simplex $\pi_0 \wedge \dots \wedge \pi_{m-1}$ and so by the Yoneda lemma there is a unique surjective map of pointed cyclic sets
\[
f_{\omega} : \Lambda^{m-1}[- ] \to \Bcy(\Pi^d, \overline{\omega})[-]
\]
mapping the generator of $\Lambda^{m-1}[-]$ to the generator $\pi_0 \wedge \dots \wedge \pi_{m-1}$. Since $\omega$ has period $s$, the generator $\pi_0 \wedge \dots \wedge \pi_{m-1}$ is fixed by the cyclic operator $t_m^s$ and so $f_\omega$ factors over the quotient subgroup of order $i = m/s$,
\[
f_\omega : (\Lambda^{m-1}[-]/C_i) \to \Bcy(\Pi^d, \overline{\omega})[-] .
\] 
From \cite[Section 7.2.]{HesselholtMadsen97}  we have a $\TT$-equivariant homeomorphism 
\[
| \Lambda^{m-1}[-]| \overset{\sim}\to \Delta^{m-1} \times \TT
\]
(where $\TT$ acts on the second factor) where the dual cyclic operator acts on $\Delta^{m-1}$ by the affine map that cyclicly permutes the vertices and on $\TT$ by rotation through $2\pi/m$. Thus $f_\omega$ gives a continuous $\TT$-equivariant surjection
\[
f_\omega : (\Delta^{m-1} \times_{C_i} \TT)_+ \to \Bcy(\Pi^d, \overline{\omega}).
\]
If $\overline{\omega}$ has no cyclic repetitions then all the faces of the generator 
\[
\pi_0 \wedge \dots \wedge \pi_{m-1}
\]
 are the basepoint $0$, and so $f_\omega$ collapses the entire boundary $\del \Delta^{m-1}$ of $ \Delta^{m-1}$ to the basepoint. There are no other relations, and so we have a $\TT$-equivariant homeomorphism
\[
f_\omega : (\Delta^{m-1} / \del \Delta^{m-1}) \wedge_{C_i} \TT_+ \overset{\sim}\to \Bcy(\Pi^d, \overline{\omega}).
\] 
Now we use the identification of the $C_m$-space $\Delta^{m-1} / \del \Delta^{m-1}$ with the one-point compactification of the reduced regular representation to finish the proof of part \emph{(1)}.

In case $\overline{\omega}$ \emph{does} have a cyclic repetition, then $\pi_0 \wedge \dots \wedge \pi_{m-1}$ will have a non-basepoint face. So $f_\omega$ collapses at least one codimension $1$ face (and its $\TT$-orbit) to the basepoint, and leaves at least one codimension $1$ face, say $F \subs \Delta^{m-1}$, un-collapsed. The cone on $F$ is canonically homeomorphic to $\Delta^{m-1}$ and has a canonical null-homotopy given by shrinking down to the basepoint of the cone. This null-homotopy then induces a null-homotopy on $\Bcy(\Pi^d, \overline{\omega})$. However, it may fail to be a $\TT$-equivariant. To get this we note that there is a $C_i$-equivariant homeomorphism
\[
\Delta^{s-1} * \cdots * \Delta^{s-1} \to \Delta^{m-1}
\]
 where $C_i$ acts on the left by cyclically permuting the factors of the join. Again $f_\omega$ will collapse at least one codimension $1$ face of $\Delta^{s-1}$ and leave at least one un-collapsed, say $F \subs \Delta^{s-1}$. Now the null-homotopy of $cone(F)$ will induce a $\TT$-equivariant null-homotopy 
 \[
 \Bcy(\Pi^d, \overline{\omega}) \wedge [0,1]_+ \to \Bcy(\Pi^d, \overline{\omega}).
 \]
This completes the proof.  
%This leads to a $\TT$-equivariant null-homotopy of $N^{cy}(\Pi^3, \overline{\omega})$ (as in \cite{Hesselholt2007} Lemma 1.6).
\end{proof}

Let $\CC(i)$ denote the $1$-dimensional complex $\TT$-representation where $z\in \TT \subs \CC$ acts through multiplication with the $i$'th power $z^i$. For $i \geq 1$ let $\lambda_i = \CC(1) \oplus \dots \oplus \CC(i)$.

\begin{lemma}\label{RegularReps}
The regular representation $\RR[C_m]$ is isomorphic to $\RR \oplus \lambda_{\frac{m-2}{2}} \oplus \RR_-$ if $m$ is even, and $\RR \oplus \lambda_{\frac{m-1}{2}}$ if $m$ is odd.
\end{lemma}
\begin{proof}
By Maschke's theorem $\RR[C_m]$ is a semisimple ring. In particular it decomposes as a sum of irreducible sub-representations. Furthermore it contains a copy of every irreducible $C_m$-representation (since for any such $V$ and any non-zero $v \in V$, there is a surjection $\RR[C_m] \to V$ given by  $\sum r_g g \mapsto \sum r_g gv$ so, by semisimplicity,  $V$ embeds into $\RR[C_m]$). Since $\RR$ and $\CC(i)$, for $1 \leq i \leq \lfloor \frac{m-1}{2} \rfloor$, (and $\RR_{-}$ in case $m$ is even) are irreducible, pair-wise non-isomorphic, and have real dimensions summing to $\dim_{\RR} \RR[C_m] = m$ this completes the proof.
%Use double coset formula, or use that $\RR[G] \cong \CC[G]^{Gal(\CC/\RR)}$.
%For any finite group $G$ we have $\RR[G] \cong \CC[G]^{\sigma}$ where $\sigma \in Gal(\CC/\RR)$ is complex conjugation. Using that
%\[
%\CC[G] \cong \bigoplus_{f \in Hom(G,S^1)} \CC_f
%\]
%[Finish: Check the isomorphism is Galois equivariant and take fixed points...]
\end{proof}

\begin{lemma}\label{FinalClassificationLemma}
Let $\overline{\omega}$ be a cyclic word with no cyclic repetitions, of length $m$, period $s$ and with $i = m/s$ blocks.
\begin{enumerate}
\item If $s$ is even, then there is a $\TT$-equivariant homeomorphism 
\[
\Sigma B(\overline{\omega}) \simeq S^{\lambda_{m/2}} \wedge (\TT/C_i)_+
\]
\item If both $s$ and $i$ are odd, then there is a $\TT$-equivariant homeomorphism 
\[
B(\overline{\omega}) \simeq S^{\lambda_{(m-1)/2}} \wedge (\TT/C_i)_+
\]
\item If $s$ is odd, and $i$ is even, then
there is a $\TT$-equivariant homeomorphism 
\[
B(\overline{\omega}) \simeq S^{\lambda_{(m-2)/2}} \wedge \RR P^2(i)
\]
where $\RR P^2(i)$ is by definition the mapping cone of the quotient map ${\TT/C_{i/2}}_+ \to {\TT/C_i}_+$.
\end{enumerate}
\end{lemma}
\begin{proof}
We use \cref{RegularReps} in each case.
\begin{enumerate}
\item If $s = 2k$ is even then so is $m$ and so 
\[
\RR[C_m]-1 \cong \lambda_{(m-2)/2} \oplus \RR_-.
\]
The restriction along the inclusion $C_i \subs C_m$ (given by $\sigma_i \mapsto \sigma_m^s$ on generators) turns the sign representation $\RR_-$ into a trivial representation. Since $\CC(\frac{m}{2}) = \CC(ki) = \RR \oplus \RR$ as $C_i$-representations, we have
\[
S^1 \wedge S^{\RR[C_m]-1} \wedge_{C_i} \TT_+ \cong  S^{\lambda_{m/2}} \wedge_{C_i} \TT_+ \cong  S^{\lambda_{m/2}} \wedge (\TT/C_i)_+
\]
where the last isomorphism uses that $\lambda_{m/2}$ extends to a representation of $\TT$ and so allows the $\TT$-equivariant untwisting map $(x,z) \mapsto (xz,zC_i)$.
\item If both $s$ and $i$ are odd then so is $m = si$ and so 
\[
\RR[C_m] - 1 = \lambda_{(m-1)/2}.
\]
Then we proceed as above, using the untwisting map.
\item If $s$ is odd and $i$ is even, then $m$ is even and so 
\[
\RR[C_m] -1 = \lambda_{(m-2)/2} \oplus \RR_ - .
\]
The restriction along the inclusion $C_i \subs C_m$ leaves the sign representation unchanged. Now repeat the argument of \cite[Cor. 7.2.]{HesselholtMadsen97}.
\end{enumerate}
\end{proof}

%%%%%%%%%%%%%%%%%%%%%%%%%%%%%%%%%%%%%%%
%%%%%%%%%%%%%%%%%%%%%%%%%%%%%%%%%%%%%%%

\subsection{Homology of $\Bcy(\Pi^d; \overline{\omega})$.}

\begin{proposition}\label{HomologyComponents}
Let $R$ be any commutative ring. Let $\overline{\omega}$ be a cyclic word with no cyclic repetitions, of length $m$, period $s$ and with $i = m/s$ blocks. The singular homology $\Htilde_*(B(\overline{\omega});R)$ is concentrated in degrees $m-1$ and $m$. If either $s$ is even, or both $s$ and $i$ are odd, then the $R$-modules in degree $m-1$ and $m$ are free of rank $1$. If $s$ is odd and $i$ is even, then the $R$-module in degree $m-1$ is isomorphic to $R/2R$, and the $R$-module in degree $m$ is isomorphic to $R[2]$, the $2$-torsion in $R$.

Furthermore, when $m$ and $s$ have the same parity, Connes' operator takes a generator $y_{\overline{\omega}}$ of the $R$-module  $\Htilde_{m-1}(B(\overline{\omega});R)$ to $i$ times a generator $z_{\overline{\omega}}$ of $\Htilde_{m}(B(\overline{\omega});R)$, that is $d(y_{\overline{\omega}}) = i z_{\overline{\omega}}$. When $s$ is odd and $i$ is even, Connes' operator acts trivially. %takes a generator $y_{\overline{\omega}} \in \Htilde_{m-1}(B(\overline{\omega};k)$ to a unit times a generator $z_{\overline{\omega}} \in \Htilde_{m}(B(\overline{\omega};k)$.
\end{proposition}
\begin{proof}
The homology computations follow directly from \cref{FinalClassificationLemma}.  We also deduce the behaviour of  Connes' operator from \cref{FinalClassificationLemma}. In general for a pointed $\TT$-space $X$, Connes' operator $d : \Htilde_*(X;k) \to \Htilde_{*+1}(X;k)$ is given by taking the cross-product with the fundamental class $[\TT]$ and then applying the map induced by the action $\mu : \TT_+ \wedge X \to X$ on reduced homology. We consider two cases.
\begin{enumerate}
\item When $X = {\TT/C_i}_+$ we claim that $d : \Htilde_0(X;k) \to \Htilde_1(X;k)$ is multiplication by $i$ (up to a unit).  By definition it is given as the composite of two maps. The first is the map
\begin{equation}\label{CrossProduct}
\Htilde_0(X) \cong \Htilde_1(\TT_+) \otimes \Htilde_0(X)
\to 
\Htilde_1(\TT_+ \wedge X) 	\qquad [z_0C_i] \mapsto [\TT \times z_0C_i]
\end{equation}
where $\TT \times z_0C_i$ is the $1$-cycle that factors through the map 
\[
\TT \cong \Delta^1/\del  
\Delta^1 \to \TT_+ \wedge X 		\qquad w \mapsto w \wedge z_0C_i
\]
Then \cref{CrossProduct} is composed with 
\[
\mu_* : \Htilde_1(\TT_+ \wedge X) \to \Htilde_1(X) .
\]
The resulting $1$-cycle in $X$ factors through the composite map
\[
\TT \to \TT_+ \wedge \TT/C_i \to \TT/C_i 	\qquad w \mapsto wz_0C_i
\]
which is a degree $i$ map.
\item When $X = \RR P^2(i)$ we claim that $d : \Htilde_1(X;k) \to \Htilde_2(X;k)$ is trivial. Consider the diagram
\begin{equation*}
\begin{tikzcd}[row sep=3em, column sep=3em]
\TT_+ \wedge (\TT/C_i)_+ \arrow[r] \arrow[d, "id \wedge q"] & (\TT/C_i)_+ \arrow[d, "q"] \\
\TT_+ \wedge \RR P^2(i) \arrow[r] & \RR P^2(i)
\end{tikzcd}
\end{equation*}
Now $H_1((\TT/C_i)_+)$ surjects onto $H_1(\RR P^2(i) )$ and the above diagram commutes. %top map factors over the surjection 
%\[H_1(\TT_+) \otimes H_1((\TT/C_i)_+) \to H_2(\TT_+ \wedge \TT/C_i)_+). 
%\]
Since $H_2((\TT/C_i)_+)=0$ the claim follows.
\end{enumerate} \end{proof}

\begin{remark}\label{RemarkForEvenPrime}
From this proposition we see that if the characteristic of $k$ is different than $2$ then $\Htilde_*(B(\overline{\omega};k)$ is trivial when $s$ is odd and $i$ is even. On the other hand, if $k$ has characteristic $2$ then $k/2k = k[2] = k$. This explains why the combinatorics of Geller, Reid, and Weibel in \cite{GellerReidWeibel}, avoids cyclic words whose period is not congruent (mod $2$) to the length \cite[Remark 3.9.1.]{GellerReidWeibel} since they work over characteristic zero fields.
\end{remark}

%%%%%%%%%%%%%%%%%%%%%%%%%%%%%%%%%%%%%%%
%%%%%%%%%%%%%%%%%%%%%%%%%%%%%%%%%%%%%%%

\subsection{$\THH$ of the coordinate axes}
We put together the various results of the previous sections to describe $\THH(A,B,I)$ as a cyclotomic spectrum.
Again, let $A = A_d, B = B_d$ and $I = I_d$. 

We will use the Segal conjecture for $C_p$. This is the statement that the map 
\[
\SSS \to \SSS^{tC_p}
\]
identifies the codomain as the $p$-completion of the domain. For a proof of this see \cite[Theorem III.1.7]{NikolausScholze}, though the result was originally proved by Lin \cite{Lin1980} (for $p=2$) and Gunawardena \cite{Gunawardena1980} (for $p$ odd) in 1980. 

\begin{lemma}\label{LaxEquivalenceLemma}
Let $T$ be a bounded below spectrum with $C_p$-action and $X$ a finite pointed $C_p$-CW-complex. Then the lax symmetric monoidal structure map 
\[
T^{tC_p} \otimes (\Sigma^\infty X)^{tC_p} \longto (T \otimes \Sigma^\infty X)^{tC_p}
\]
is an equivalence.
\end{lemma}
\begin{proof}
Since both $T^{tC_p} \otimes (-)^{tC_p}$ and $(T \otimes -)^{tC_p}$ are exact functors we may reduce to checking the statement for the $C_p$-spectra $\SSS$ and $\SSS \otimes {C_p}_+$. This is because  $\Sigma^\infty X$ may be constructed from $\SSS$ and $\SSS \otimes {C_p}_+$ using finitely many cofiber sequences (since $X$ is built by attaching finitely many $C_p$-cells). Replacing $\Sigma^\infty X$ by $\SSS$ the map in question reduces to
\[
T^{tC_p} \otimes \SSS^{\widehat{}}_p \longto T^{tC_p}
\]
where we use the Segal conjecture to identify $\SSS^{tC_p} \simeq \SSS^{\hat{}}_p$. Since $T$ is bounded below it follows from \cite[Lemma I.2.9]{NikolausScholze}) that $T^{tC_p}$ is $p$-complete and so the map is an equivalence.  Replacing $\Sigma^\infty X$ by $\SSS \otimes {C_p}_+$ instead we see that both the domain and codomain of the map are zero, since $(-)^{tC_p}$ kills induced spectra as well as spectra of the form $T \otimes Z$ where $Z$ is induced (cf. \cite[Lemma I.3.8. (i) and (ii)]{NikolausScholze}). 
\end{proof}

We will describe the Frobenius map on 
\[
\THH(k(\Pi^d)) \simeq \THH(k) \otimes \Bcy(\Pi^d)
\]
in terms of the splitting of $\Bcy(\Pi^d)$ into the $\TT$-equivariant subspaces $B(\overline{\omega})$. The unstable Frobenius $\psi_p : \Bcy(\Pi^d) \to \Bcy(\Pi^d)^{C_p}$ (defined in \cref{UnstableFrobenius}) restricts to a homeomorphism 
\[
\psi_p : B(\overline{\omega}) \to B(\overline{\omega^{\star p}})^{C_p}
\]
landing in the subspace $B(\overline{\omega^{\star p}})$ corresponding to the cyclic word $\overline{\omega^{\star p}}$ which has length $pm$ and period $s$ (if $\overline{\omega}$ has length $m$ and period $s$).

\begin{proposition}\label{CyclotomicStructure}
There is a $\TT$-equivariant equivalence of spectra
\[
\THH(A,B,I) \simeq \bigoplus \THH(k) \otimes B(\overline{\omega})
\]
where the sum is indexed over cyclic words of length $m \geq 2$, having no cyclic repetitions. 
Under this equivalence the Frobenius map restricts to the map
\[
\THH(k) \otimes B(\overline{\omega}) \overset{\phi_p \otimes \tilde{\phi}_p}\longto \THH(k)^{tC_p} \otimes B(\overline{\omega^{\star p}})^{tC_p} \longto ( \THH(k) \otimes B(\overline{\omega^{\star p}}) )^{tC_p}
\]
where the second map is the lax symmetric monoidal structure of the Tate-$C_p$-construction. This second map is an equivalence, while the restricted Frobenius $\tilde{{\phi}_p} : \Sigma^\infty B(\overline{\omega}) \to (\Sigma^\infty B(\overline{\omega^{\star p}}))^{tC_p}$ is a $p$-adic equivalence.
\end{proposition}
%\begin{prop}\label{CyclotomicStructure}
%There is a $\TT$-equivariant equivalence of spectra
%\[
%\THH(A,B,I) \simeq \bigoplus \THH(k) \otimes B(\overline{\omega})
%\]
%where the sum is indexed over cyclic words of length $m \geq 2$, having no cyclic repetitions. 
%Under this equivalence the Frobenius map restricts to the map
%\begin{eqnarray*}
%\THH(k) \otimes B(\overline{\omega} ) 
%&\overset{\phi \otimes \psi_p}{\longto} \THH(k)^{tC_p} \otimes B(\overline{\omega^{\star p}})^{C_p}
%\\ &\overset{\simeq}\longto (\THH(k) \otimes B( \overline{\omega^{\star p}})^{C_p})^{tC_p} 
%\\ &\overset{\simeq}\longto \left( \THH(k) \otimes B(\overline{\omega^{\star p}})  \right)^{tC_p}  
%\end{eqnarray*}%\todo{identify the Frobenius}
%where $\phi : \THH(k) \to \THH(k)^{tC_p}$ is the canonical Frobenius on $\THH(k)$. The second map moves the trivial $C_p$-space $B( \overline{\omega} \star \dots \star \overline{\omega})^{C_p}$ inside the Tate $C_p$-construction. The final map is induced by the inclusion of $B( \overline{\omega} \star \dots \star \overline{\omega})^{C_p}$ into $B( \overline{\omega} \star \dots \star \overline{\omega})$. 
%\end{prop}

%We abuse notation and continue to denote the composite in \cref{CyclotomicStructure} by $\phi$. 
\begin{proof}
Applying \cref{LaxEquivalenceLemma} with $T = \THH(k)$ and $X = B(\overline{\omega^{\star p}})$ we get that the map
\[
\THH(k)^{tC_p} \otimes B(\overline{\omega^{\star p}})^{tC_p} \longto ( \THH(k) \otimes B(\overline{\omega^{\star p}}) )^{tC_p}
\]
is an equivalence.
It remains to show that 
\[
\tilde{\phi}_p : \SSS \otimes B(\overline{\omega}) \to (\SSS \otimes B(\overline{\omega^{\star p}}))^{tC_p}\]
is a $p$-adic equivalence. To do this we factor it as follows. To ease notation, let $Y = |B(\overline{\omega^{\star p}})[-]|$. By definition $\tilde \phi_p$ factors as
\[
\SSS \otimes B(\overline{\omega}) \overset{\tilde{\Delta_p}}\longto \SSS \otimes |\sd_p B(\overline{\omega^{\star p}})[-]|^{C_p} \overset{D_p}\longto \SSS \otimes Y^{C_p} \overset{\tilde{\phi'_p}}\longto  (\SSS \otimes Y)^{tC_p}
\]
where $\tilde{\Delta_p}$ (the space-level diagonal) and $D_p$  are homeomorphisms as remarked in \cref{UnstableFrobenius}. The final map (denoted by $\tilde{\phi'_p}$)  is the composition \[
\tilde{\phi'_p} : \SSS \otimes Y^{C_p} \overset{\gamma}\longto \SSS \otimes Y^{hC_p} \longto (\SSS \otimes Y)^{hC_p} \overset{\can}\longto (\SSS \otimes Y)^{tC_p}
\]
This map fits into the following commutative diagram
\begin{center}
\begin{tikzcd}[row sep=4em, column sep=3em]
\SSS \otimes Y^{C_p} \arrow[dd, swap, "\Delta_p \otimes id"] \arrow[rr, "\tilde{\phi'_p}"] \arrow[dr] &  & (\SSS \otimes Y)^{tC_p}  \\
\ & (\SSS \otimes Y^{C_p})^{hC_p} \arrow[d, "\can"] \arrow[ur , "\can \circ\ i"]  & \  \\
\SSS^{tC_p} \otimes Y^{C_p} \arrow[r, "(1.)"] & (\SSS \otimes Y^{C_p})^{tC_p} \arrow[uur, bend right=35, swap, "(2.)"] &  \
\end{tikzcd}
\end{center}
Here $\Delta_p : \SSS \to \SSS^{tC_p}$ is a $p$-adic equivalence by the Segal conjecture.  The map labelled $(1.)$ is the equivalence witnessing that $(-)^{tC_p}$ is an exact functor, and $B(\overline{\omega^{\star p}})^{C_p}$ is finite and has trivial $C_p$-action.  The left-most square commutes by construction of the map $(1.)$. Indeed by exactness in the variable $Y^{C_p}$ (a spectrum with trivial $C_p$-action) one reduces to the case $Y^{C_p} = \SSS$ where the square becomes
\begin{center}
\begin{tikzcd}[row sep=3em, column sep=3em]
\SSS \arrow[d, "\Delta_p"] \arrow[r] & \SSS^{hC_p} \arrow[d, "\can"] \\
\SSS^{tC_p} \arrow[r, " \id "] & \SSS^{tC_p}
\end{tikzcd}
\end{center}

 The map labelled $(2.)$ is the equivalence witnessing that the inclusion of the $C_p$-singular set $B^{C_p} \subs B$ induces an equivalence 
\[
(T \otimes B^{C_p})^{tC_p} \to (T \otimes B)^{tC_p}
\]
 for any $C_p$-spectrum $T$ and finite $C_p$-CW-complex $B$, cf. \cite[Lemma 9.1.]{HesselholtMadsen97}. The right-most triangle (with $(2.)$ as a side) commutes since $\can$ is natural. 
 
Finally, in the top triangle, the map $\SSS \otimes Y^{C_p} \to (\SSS \otimes Y^{C_p})^{hC_p}$ arises since $\SSS \otimes Y^{C_p}$ has trivial $C_p$-action. The homotopy $C_p$-fixed point functor is right adjoint to the functor equipping a spectrum with a trivial $C_p$ action. In particular, the map $(\SSS \otimes Y)^{hC_p} \to \SSS \otimes Y$ is the universal map from a spectrum with trivial $C_p$-action to the target. It follows that the  following diagram commutes.
\begin{center}
\begin{tikzcd}[row sep=3em, column sep=3em]
\SSS \otimes Y^{C_p} \arrow[dr] \arrow[r] & \SSS \otimes Y^{hC_p} \arrow[r] & (\SSS \otimes Y)^{hC_p} \\
\ & ( \SSS \otimes Y^{C_p})^{hC_p} \arrow[ur] & \
\end{tikzcd}
\end{center}
This completes the proof.
\end{proof}

I thank Malte Leip for abundant help with the above argument.

\begin{corollary}\label{FrobeniusConnectivity}
Let $\overline{\omega}$ be a cyclic word having no cyclic repetitions and length $m$. The map \[
\phi_p : \THH(k) \otimes B(\overline{\omega}) \to (\THH(k) \otimes B( \overline{\omega^{\star p}}) )^{tC_p}
\]
induces an isomorphism on homotopy groups in degrees greater than or equal to $m$ %\todo[line]{check that this is true, should this rather be stated in terms of homotopy groups?}
\end{corollary}
\begin{proof}
The Frobenius map $\phi_p : \THH(k) \to \THH(k)^{tC_p}$ induces an isomorphism on non-negative homotopy groups. This is shown in \cite[Prop. IV. 4.13]{NikolausScholze} for $k = \FF_p$, and in \cite[Section 5.5]{HesselholtMadsen97} for $k$ a perfect field of characteristic $p$. %By \cite[Section 2]{BHM} $\psi_p$ is an equivalence. 
From \cref{CyclotomicStructure} we know that 
\[
\phi_p : \THH(k) \otimes B(\overline{\omega}) \to (\THH(k) \otimes B( \overline{\omega^{\star p}}) )^{tC_p}
\]
factors as the composite of an equivalence and the tensor product of the Frobenius  on $\THH(k)$,  $\phi_p : \THH(k) \to \THH(k)^{tC_p}$ with the map $\tilde{{\phi}_p} : \Sigma^\infty B(\overline{\omega}) \to (\Sigma^\infty B(\overline{\omega^{\star p}}))^{tC_p}$ which is a $p$-adic equivalence. 
The result now follows from  \cref{HomologyComponents}  as can be shown for example with the Atiyah-Hirzebruch spectral sequence. Indeed, given a map of spectra $\phi: T \to T'$ which is an equivalence in high degrees, and a space $B$ with homology bounded above, then $\phi \otimes B: T \otimes B \to T' \otimes B$ is an equivalence in high degrees. To see this, take the cofiber $\mathrm{cof}(\phi \otimes B) = \mathrm{cof}(\phi) \otimes B$, we must show it is highly coconnective. The Atiyah-Hirzebruch spectral sequence for $\pi_*(\mathrm{cof}(\phi) \otimes B)$ has $E^2$-page 
\[
E_{i,j} = \Htilde_i(B; \pi_j(\mathrm{cof}(\phi))).
\]
When $i$ is large $E_{i,j} = 0$ and when $j$ is large $E_{i,j} = 0$ so there is some diagonal beyond which all groups are zero. In the case at hand this diagonal is determined by \cref{HomologyComponents} and is given by $i+j = m$.
\end{proof}

\noindent By \cref{HomologyComponents} the homology of $B(\overline{\omega})$ depends only on the length and period of $\overline{\omega}$. For $\overline{\omega}$ of length $m$ and period $s$ we therefore introduce the notation $B(m,s) = B(\overline{\omega})$. The function $\cyc_d(s)$ counts how many cyclic words with no cyclic repetitions, of length $s$, and period $s$, there are. For a fixed $m$ and $s$ a divisor of $m$, there are $\cyc_d(s)$ many cyclic words with no cyclic repetitions, of length $m$, and period $s$. Since the homology of $B(m,s)$ depends only on the parity of $m$ and $s$ (by \cref{HomologyComponents}) we have the following equivalences; if $p > 2$ then
\begin{eqnarray} \label{BiRelativeTHH}
\THH(A,B,I) \simeq \bigoplus_{\substack{m \geq 2 \\ \text{even}}} \bigoplus_{\substack{s \mid m \\ \text{even}}} (\THH(k) \otimes B(m,s))^{\oplus \cyc_d(s)}
\\ \oplus
\bigoplus_{\substack{m \geq 2 \\ \text{odd}}} \bigoplus_{\substack{s \mid m \\ \text{odd}}} (\THH(k) \otimes B(m,s))^{\oplus \cyc_d(s)} .
\end{eqnarray}
If $p=2$ then we add the similar double sum indexed over $m$ even and $s$ odd, i.e.\ 
\begin{eqnarray*} \label{BiRelativeTHHp=2}
\THH(A,B,I) \simeq \bigoplus_{\substack{m \geq 2 \\ \text{even}}} \bigoplus_{\substack{s \mid m \\ \text{even}}} (\THH(k) \otimes B(m,s))^{\oplus \cyc_d(s)}
\\ \oplus
\bigoplus_{\substack{m \geq 2 \\ \text{odd}}} \bigoplus_{\substack{s \mid m \\ \text{odd}}} (\THH(k) \otimes B(m,s))^{\oplus \cyc_d(s)}
\\ \oplus \bigoplus_{\substack{m \geq 2 \\ \text{even}}} \bigoplus_{\substack{s \mid m \\ \text{odd}}} (\THH(k) \otimes B(m,s))^{\oplus \cyc_d(s)} .
\end{eqnarray*}

Note that we do not need to know the homotopy-type of $B(m,s)$ for the above equivalences, since the homotopy type of $\THH(k) \otimes B(m,s)$ is determined from the homology of $B(m,s)$. Indeed the Atiyah-Hirzebruch spectral sequence for the $\THH(k)$-homology of $B(m,s)$ degenerates at the $E_2$-page with trivial extension problems, and thus computes $\pi_*(\THH(k)\otimes B(m,s))$. It is the simplicity of the homology of $B(m,s)$ that makes it possible for us to compute this Atiyah-Hirzebruch spectral sequence.

%%%%%%%%%%%%%%%%%%%%%%%%%%%%%%%%%%%%%%%%
%%%%%%%%%%%%%%%%%%%%%%%%%%%%%%%%%%%%%%%%
%%%%%%%%%%%%%%%%%%%%%%%%%%%%%%%%%%%%%%%%

%%%%%%%%%%%%%%%%%%%%%%%%%%%%%%%%%%%%%%%%
%%%%%%%%%%%%%%%%%%%%%%%%%%%%%%%%%%%%%%%%
%%%%%%%%%%%%%%%%%%%%%%%%%%%%%%%%%%%%%%%%

\section{Negative - and periodic topological cyclic homology}\label{TPsection}
\noindent
The results of the previous section give us a good understanding of $\THH(A,B,I)$ as a cyclotomic spectrum, i.e.\ as a spectrum with $\TT$-action and with the $\TT$-equivariant Frobenius map $\phi_p : THH(A,B,I) \to \THH(A,B,I)^{tC_p}$. In this section we take the next step towards computing $\TC(A,B,I)$ and hence towards \cref{TCCoordinateAxes}. To do this we will use the results of \cite{NikolausScholze} which give a formula for $\TC$ in terms of the cyclotomic structure on $\THH$. Let us briefly recall this formula. Given a ring $R$ we view $\THH(R)$ as a cyclotomic spectrum and form the topological negative cyclic homology,
\[
\TC^-(R) = \THH(R)^{h\TT}
\]
and topological periodic cyclic homology
\[
\TP(R) = \THH(R)^{t\TT}.
\]
There is a canonical map from the former to the latter
\[
\can : \TC^-(R) \to \TP(R).
\]
Assuming $R$ is characteristic $p > 0$ it follows from an application of the Tate orbit lemma \cite[Lemma I.2.1]{NikolausScholze} that $\phi_p$ induces a map, which, by abuse of notation, we still call $\phi_p$,
\[
\phi_p : \TC^-(R) \to \TP(R).
\]
Finally, by \cite[Section II.1]{NikolausScholze} there is a cofiber sequence
\begin{equation}\label{TCcofiberSequence}
\TC(R;p) \to \TC^-(R) \overset{\phi_p - \can}\longto \TP(R).
\end{equation}
All of this readily works in the birelative case so that we also have a cofiber sequence for $\TC(A,B,I)$. 
The goal of this section is thus to compute $\TC^-(A,B,I)$ and $\TP(A,B,I)$. 

We will need the following general lemma about the Tate $\TT$-construction.

\begin{lemma}\label{TPproducts}
Let $\{X_i\}_{i\geq 0}$ be a sequence of spectra with a $\TT$-action, such that the connectivity of (the underlying spectrum) $X_i$ is unbounded, as $i$ grows. Then the canonical map $(\bigoplus_i X_i)^{t\TT} \to \prod (X_i)^{t\TT}$ is an equivalence.
\end{lemma}
\begin{proof}
Because of the connectivity assumption the canonical map $\bigoplus X_i \to \prod X_i$ is an equivalence. Consider the norm cofiber sequence
\[
(\Sigma X)_{h\TT} \to X^{h\TT} \to X^{t\TT}
\]
The first term commutes with colimits, the second with limits, so we get the cofiber sequence
\[
\bigoplus (\Sigma X_i)_{h\TT} \to \prod X_i^{h\TT} \to (\bigoplus X_i)^{t\TT}
\]
Now $( - )_{hG}$ preserves connectivity for any group $G$. To see this one may use the homotopy orbit spectral sequence whose $E^2$-page consists of ordinary group homology. Thus the cofiber sequence becomes
\[
\prod (\Sigma X_i)_{h\TT} \to \prod X_i^{h\TT} \to (\bigoplus X_i)^{t\TT}
\]
and we see that $(\bigoplus_i X_i)^{t\TT} \to \prod (X_i)^{t\TT}$ is an equivalence.
\end{proof}

%%%%%%%%%%%%%%%%%%%%%%%%%%%%%%%%%%%%%%%%%%%%%%
%%%%%%%%%%%%%%%%%%%%%%%%%%%%%%%%%%%%%%%%%%%%%%

\subsection{The Tate spectral sequence} \label{TateSS}

Let $X$ be a spectrum with $\TT$-action. The Tate construction $\TP(X) = X^{t\TT}$ is the target of the Tate spectral sequence for the circle group $\TT$, see \cite[end of section 4.4]{HesselholtMadsen2003}, and also \cite[section 3]{BlumbergMandell2017}.  This spectral sequence has the form
\[
E^2(\TT,X) = S\{t^{\pm 1}\} \otimes \pi_*X \quad \Rightarrow \quad \pi_* \TP(X)
\]
Here $t$ has bi-degree $(-2,0)$ and is a generator  of $H^2(\CC P^\infty ; \ZZ)$. %(i.e.\ the first Chern class of $\mathcal{O}$.)
This spectral sequence is conditionally convergent. When $X = \THH(R)$ is an $E_\infty$ algebra the spectral sequence is multiplicative, and if $R$ is a $k$-algebra then $E^r(\TT,\THH(R))$ is a module over $E^r(\TT,\THH(k))$. By B\"okstedt periodicity, when $k$ is a perfect field of characteristic $p>0$ there is an isomorphism 
\[
\THH_*(k) \simeq k[x]
\]
with $x$ in degree $2$. In particular, in this case $\TP_*(R)$ is periodic. 
The following is well-known.

\begin{proposition}\label{TPperiodicity}
Let $k$ be a perfect field of characteristic $p > 0$, then there is an isomorphism
\[ 
\TP_*(k) \simeq W(k)[t^{\pm 1}]
\]
\end{proposition}
\begin{proof}
By B\"okstedt periodicity the $E^2$-term of the spectral sequence takes the form $k[t^{\pm 1}] \otimes k[x]$, and since all classes are in even total degree there are no non-trivial differentials on any page. Thus $E^2 = E^\infty$ and it remains to determine the extensions. Again by periodicity it suffices to show that $\TP_0(k) = W(k)$. For $k = \FF_p$ this is done in \cite[Cor. IV. 4.8.]{NikolausScholze}. To conclude the result for $k$ we use functoriality along the map $\FF_p \to k$.  More precisely we argue as follows. From the Tate spectral sequence we have a complete descending multiplicative filtration
\[
\dots \subs \Fil^{i+1}(k) \subs \Fil^i(k) \subs \dots \Fil^1(k) \subs \Fil^0 = \TP_0(k)
\]
with  associated graded $\gr^i(k) \simeq \THH_{2i}(k) \simeq k$. By the universal property of the $p$-typical Witt vectors $W(k)$ \cite[Chapter II, paragraph 5]{SerreLocalFields} we get a unique multiplicative continous map $W(k) \to \TP_0(k)$. By functoriality we have a commuting diagram
\begin{equation*}
\begin{tikzcd}[row sep=3em, column sep=3em]
W(\FF_p) \arrow[r, "\simeq"] \arrow[d] & \TP_0(\FF_p) \arrow[d] \\
W(k) \arrow[r] & \TP_0(k)
\end{tikzcd}
\end{equation*}
of ring homomorphisms. In particular the map on associated graded induced by $W(k) \to \TP_0(k)$ maps $1$ to $1$ and so must be an isomorphism. It follows that $W(k) \to \TP_0(k)$ is itself an isomorphism.
\end{proof}

%\begin{sseqdata}[ name = TateSS,% xscale = 1, yscale = 1, 
%homological Serre grading , classes = { draw = none } ]
%\class["t"](-2,0) \class(0,0) \class[](0,2) \class(-2,1)
%\d2(0,0) 
%\end{sseqdata}
%\printpage[name = TateSS, page = 2] 

\begin{lemma}\label{TateSpectralSequenceProperties}
Let $k$ be a perfect field of characteristic $p > 0$. The elements  $x$ and $t^{\pm}$ from the $E^2$-page of the Tate spectral sequence are infinite cycles.
\end{lemma}
\begin{proof}
This is true for degree reasons, by B\"okstedt periodicity. 
\end{proof}

\begin{lemma}\label{SecondTateDifferential}
Let $X$ be a spectrum with $\TT$-action such that the underlying spectrum is an $H\ZZ$-module . The $d^2$ differential of the $\TT$-Tate spectral sequence is given by $d^2(\alpha) = t d(\alpha)$ where $d$ is Connes' operator.
\end{lemma}
\begin{proof}
See \cite[Lemma 1.4.2]{Hesselholt96}  %[or B\"oksted-Madsen section 5?]
\end{proof}

\noindent Following \cref{HomologyComponents} we now choose some generators for the homology of the spaces $B(m,s)$. When $m$ and $s$ have the same parity let $z_{(m,s)}$ be a generator for $\Htilde_m(B(m,s);k)$ and let  $y_{(m,s)}$ a generator for $\Htilde_{m-1}(B(m,s);k)$. When $m$ is even and $s$ is odd  let $z_{(m,s)}$ be a generator of $\Htilde_{m-1}(B(m,s);k)$. Finally, when $p=2$, so $\Htilde_{m}(B(m,s);k)$ is free of rank $1$ over $k$,  let $w_{(m,s)}$ be a generator in degree $m$.

\begin{lemma}\label{InfiniteCycleLemma}
As an element of the $E^2$-page of the Tate spectral sequence, $z_{(m,s)}$ is an infinite cycle, i.e.\ it is annihilated by all differentials.
\end{lemma}
\begin{proof}
By \cite[Theorem B]{HesselholtMadsen97} for any perfect field $k$ of positive characteristic, there is an equivalence $\tau_{\geq 0} \TC(k) \simeq \ZZ_p$. As a result we obtain a map 
\[
\ZZ_p \simeq \tau_{\geq 0} \TC(k) \to \TC(k) \to \TC^{-}(k) \to \THH(k) .
\]
of spectra. The map $\TC(k) \to \TC^-(k)$ is the first map of \cref{TCcofiberSequence}. The second map is the canonical map $(-)^{h\TT} \to (-)$. If we equip the all terms except $\THH(k)$ with trivial $\TT$-action then each map is $\TT$-equivariant, and so the composition yields a $\TT$-equivariant map 
\[
\ZZ_p \to \THH(k).
\]
Since it is $\TT$-equivariant, this map induces a map on Tate spectral sequences, which when $p$ is odd takes the following form,
\begin{equation*}
\begin{tikzcd}[row sep=3em, column sep=3em]
E^2 = \ZZ_p[t^{\pm 1}] \{ y_{(m,s)} , z_{(m,s)}\} \arrow[r, Rightarrow] \arrow[d, ] & \pi_*(\ZZ_p \otimes B(m,s))^{t\TT} \arrow[d , ] \\
E^2 = k[t^{\pm 1}, x]  \{ y_{(m,s)} , z_{(m,s)}\} \arrow[r, Rightarrow] & \pi_*(\THH(k) \otimes B(m,s))^{t\TT}
\end{tikzcd}
\end{equation*}
For degree reasons $z_{(m,s)}$ is an infinite cycle in the top spectral sequence. It follows that $z_{(m,s)}$, viewed as a class in the bottom spectral sequence, is an infinite cycle. When $m$ is even and $s$ is odd, and $p=2$, the spectral sequences in question have a slightly different form, with generators $z_{(m,s)}$ and $w_{(m,s)}$ (for the $\THH(k)$-homology, cf.\ notation introduced in the paragraph before the lemma statement). However, the argument still works, i.e.\  for degree reasons $z_{(m,s)}$ must be an infinite cycle in the top spectral sequence, which forces the same conclusion for the target spectral sequence.
%Using \cref{SecondTateDifferential} and \cref{HomologyComponents}  we see that the $E^3$ page of the top spectral sequence is $\ZZ_p/i\ZZ_p[t^{\pm 1}]\{ z_{(m,s)} \}$, hence has no more differentials. It follows that $z_{(m,s)}$, viewed as a class in the bottom spectral sequence, is an infinite cycle.
\end{proof}

%\begin{lemma} % This is dealt with in the above calculation, TP(k) = W(k).
%The extensions in the Tate spectral sequence for THH-stuff are maximal... cf. Hesselholt-Madsen 1997, and Hesselholt notes at MIT...
%\end{lemma}

\subsection{$\TC^-$ and $\TP$ of coordinate axes}

In this section we compute $\TC^-$ and $\TP$ using the Tate spectral sequence.
First assume $p > 2$. By \cref{TPproducts} and \cref{BiRelativeTHH} we have
\begin{eqnarray*}
\TP(A,B,I) \simeq \prod_{\substack{m \geq 2 \\ \text{even}}} \prod_{\substack{s \mid m \\ \text{even}}} \left( (\THH(k) \otimes B(m,s))^{t\TT} \right)^{\oplus \cyc_d(s)}
\\ \oplus
\prod_{\substack{m \geq 2 \\ \text{odd}}} \prod_{\substack{s \mid m \\ \text{odd}}} \left( (\THH(k) \otimes B(m,s))^{t\TT} \right)^{\oplus \cyc_d(s)}
\end{eqnarray*}

and likewise for negative topological cyclic homology, we have
\begin{eqnarray*}
\TC^-(A,B,I) \simeq \prod_{\substack{m \geq 2 \\ \text{even}}} \prod_{\substack{s \mid m \\ \text{even}}} \left( (\THH(k) \otimes B(m,s))^{h\TT} \right)^{\oplus \cyc_d(s)}
\\ \oplus
\prod_{\substack{m \geq 2 \\ \text{odd}}} \prod_{\substack{s \mid m \\ \text{odd}}} \left( (\THH(k) \otimes B(m,s))^{h\TT} \right)^{\oplus \cyc_d(s)}
\end{eqnarray*}

When $p=2$ there is an extra double product indexed over $m \geq 2$ even and $s \mid m$ odd (cf. \cref{HomologyComponents} and \cref{RemarkForEvenPrime}). We handle this case  in \cref{Casep=2} below.

Since $\THH(k)$ is $p$-complete, it follows from \cite[Lemma II. 4.2.]{NikolausScholze} that we may identify the homotopy $\TT$-fixed points of the Frobenius morphism for $\THH(A,B,I)$ with the map induced by the product of the maps
\[
\phi_p(m,s) : (\THH(k) \otimes B(m,s))^{h\TT}  \to (\THH(k) \otimes B(pm,s))^{t\TT}  .
\]
 Since the homotopy fixed point functor preserves co-connectivity, it follows from \cref{FrobeniusConnectivity} that $\pi_*\phi_p(m,s)$ is an isomorphism when $* \geq m$. Indeed, the fiber of 
 \[
 (\THH(k) \otimes B(m,s))  \to (\THH(k) \otimes B(pm,s))^{tC_p}
 \]
has no homotopy groups above degree $m-1$, hence the same is true of the $\TT$-homotopy fixed point spectrum.
 
\begin{proposition}\label{TPandTCminus} Let $k$ be a perfect field of characteristic $p > 0$.
Let $m \geq 2$ and $s \mid m$. Write $m = p^vm'$ and $s = p^u s'$ with $m'$ and $s'$ coprime to $p$. 
\begin{enumerate}
\item If  both $m$ and $s$ are even then $\pi_*((\THH(k) \otimes B(m,s))^{t\TT})$ as well as $\pi_*((\THH(k) \otimes B(m,s))^{h\TT})$ are concentrated in even degrees and given by
\[
\pi_{2r}((\THH(k) \otimes B(m,s))^{t\TT}) \simeq W_{v-u}(k) 
\]
and
\[
 \pi_{2r}((\THH(k) \otimes B(m,s))^{h\TT}) \simeq \left\{
        \begin{array}{ll}
           W_{v-u+1}(k)  & 2r \geq m   \\
            W_{v-u}(k) & \quad  2r < m
        \end{array}
    \right. .
\]
\item If both $m$ and $s$ are odd, then $\pi_*((\THH(k) \otimes B(m,s))^{t\TT})$ as well as $\pi_*((\THH(k) \otimes B(m,s))^{h\TT})$ are concentrated in odd degrees and given by
\[
\pi_{2r+1}((\THH(k) \otimes B(m,s))^{t\TT}) \simeq W_{v-u}(k) 
\]
and
\[
 \pi_{2r+1}((\THH(k) \otimes B(m,s))^{h\TT}) \simeq \left\{
        \begin{array}{ll}
           W_{v-u+1}(k)  &  2r \geq m  \\
            W_{v-u}(k) & \quad  2r < m
        \end{array}
    \right. .
\]
\end{enumerate}
\end{proposition}
\begin{proof}
%We apply induction on $v \geq 0$ using the Tate spectral sequence. In general the Tate spectral sequence converging to $\pi_*((\THH(k) \otimes B(m,s))^{t\TT})$ has $E^2$-page
Suppose first that $m$ and $s$ are even.
We proceed by induction on $v \geq 0$. Suppose $v = 0$, so $m = m'$ and $s = s'$. Consider the Tate spectral sequence (cf. \cref{TateSS} )
\[
E^2 = k[t^{\pm 1}, x]\{y_{(m',s')} , z_{(m',s')} \} \quad \Rightarrow \quad \pi_*(\THH(k) \otimes B(m',s') )^{t\TT}
\]
By \cref{InfiniteCycleLemma} and \cref{TateSpectralSequenceProperties} the differential structure is determined by the differentials on $y_{(m',s')}$. Furthermore 
\[
d^2(y_{(m',s')}) \overset{.}= t d(y_{(m',s')}) \overset{.}= t i z_{(m',s')}
\]
 by \cref{SecondTateDifferential} and \cref{HomologyComponents}, where $i = m'/s'$. Since $i$ is a unit in $k$, $d^2$ is an isomorphism. The figure below displays the $E^2$-page (where we have dropped the indices for clarity). All the arrows indicate isomorphisms.

\begin{figure}[b]
\centering
    \begin{tikzpicture}[scale = 0.7, every node/.style={transform shape}]
    
    \draw (-7.5,-0.5) -- (6.8,-0.5);
    \draw (0.4,-0.5) -- (0.4,5);
    \node at (0,0) {$y$};
    \node at (0,1) {$z$};
    \node at (0,2) {$xy$};
     \node at (0,3) {$xz$};
      \node at (0,4) {$\vdots$};
    %\node at (-1,0) {$u_1$};
    \node at (-2,0) {$ty$};
    \node at (2,0) {$t^{-1}y$};
     \node at (4,0) {$t^{-2}y$};
    \node at (5,0) {$\cdots$};
    %\node at (-3,0) {$tu_1$};
    \node at (-3.5,0) {$\cdots$};
    %\node at (1,0) {$t^{-1}u_1$};
   % \node at (1.5,0) {$\cdots$};

    \draw[->] (-2.3,0.2) -> (-3.8,0.9);
    \draw[->] (-0.3,0.2) -> (-1.8,0.9);
    \draw[->] (1.7,0.2) -> (0.2,0.9);
    \draw[->] (3.7,0.2) -> (2.2,0.9);    
       \draw[->] (-2.3,2.2) -> (-3.8,2.9);
    \draw[->] (-0.3,2.2) -> (-1.8,2.9);
    \draw[->] (1.7,2.2) -> (0.2,2.9);
    \draw[->] (3.7,2.2) -> (2.2,2.9);    
    \node at (-4,4) {$\cdots$};
    \node at (1.5,4) {$\cdots$};
    \node at (-7,4) {$\cdots$};
    \end{tikzpicture}
\end{figure}

Thus $E^3$, hence $E^\infty$, is trivial as claimed. To determine the $\TT$-homotopy fixed points, we truncate the Tate spectral sequence, removing the first quadrant. Thus the class $z_{(m',s')}$ and its multiples by $x^n$, are no longer hit by differentials and so 
\[
E^3 = E^\infty = k[x]\{z_{(m',s')}\}
\]
 where $z_{(m',s')}$ has degree $m$. This proves the claim for $v = 0$. The same argument works for $v > 0$ and $u = v$, since in this case $i$ is again coprime to $p$.

Suppose the claim is known for all integers less than or equal to $v$ and all $u \leq v$. As we saw above, the homotopy fixed points of the Frobenius map
\[
\pi_* (\THH(k) \otimes B(p^vm',p^us') )^{h\TT} \to \pi_* (\THH(k) \otimes B(p^{v+1}m',p^us'))^{t\TT}
\]
is an isomorphism when $* \geq p^vm'$. The induction hypothesis implies that the domain is isomorphic to $W_{v-u+1}(k)$ when $* = 2r \geq p^vm'$.   Since the Tate construction admits the structure of a lax symmetric monoidal functor \cite[Corollary I.4.3.]{NikolausScholze} it follows that $(\THH(k) \otimes B(p^{v+1}m', p^us'))^{t\TT}$ inherits a module structure over $\THH(k)^{t\TT}$. In particular it follows from the calculation 
\[
\pi_*(\THH(k)^{t\TT}) \simeq W(k)[t^\pm]
\]
(\cref{TPperiodicity}) that $\pi_*(\THH(k) \otimes B(p^{v+1}m', p^us'))^{t\TT}$ is periodic as well. Using this periodicity we conclude that
\[
\pi_{2r} (\THH(k) \otimes B(p^{v+1}m', p^us'))^{t\TT} \simeq W_{v-u+1}(k)
\]
for any $r \in \ZZ$. Considering again the Tate spectral sequence we see that we must have
\[
d^{2(v+1-u) + 2}(y_{(p^{v+1}m',p^us')}) \overset{.}= t z_{(p^{v+1}m',p^us')} (xt)^{v+1-u}
\]
(see the figure below, which is page $E^8$ when $v=2$ and $u=0$) and so we conclude that
$E^{2(v+1-u) + 3} = E^\infty$.

%\begin{comment}
\begin{figure}[b]
\centering
    \begin{tikzpicture}[scale = 0.70, every node/.style={transform shape}]
    
    \draw (-7.5,-0.5) -- (6.8,-0.5);
    \draw (0.4,-0.5) -- (0.4,5);
    \node at (0,0) {$y$};
    \node at (0,1) {$z$};
    \node at (0,2) {$xy$};
    \node at (0,3) {$xz$};
    \node at (0,4) {$x^2y$};
    \node at (0,5) {$\vdots$};
    %\node at (-1,0) {$u_1$};
    \node at (-2,0) {$ty$};
    \node at (4,0) {$t^{-2}y$};
    \node at (2,0) {$t^{-1}y$};
    \node at (5,0) {$\cdots$};
    %\node at (-3,0) {$tu_1$};
    \node at (-4.5,0) {$\cdots$};
    \node at (-10,9) { };
    %\node at (1,0) {$t^{-1}u_1$};
   % \node at (1.5,0) {$\cdots$};

    \draw[->] (-2.3,0.2) -> (-9.8,6.9); % differential on t
    \draw[->] (-0.3,0.2) -> (-7.8,6.9);  % differential on y
    \draw[->] (1.7,0.2) -> (-5.8,6.9); % differential on t^-1y
    \draw[->] (3.7,0.2) -> (-3.8,6.9);    % differential on t^-2y
    \draw[->] (-2.3,2.2) -> (-9.8,8.9); % differential on txy^2
    \draw[->] (-0.3,2.2) -> (-7.8,8.9); % differential on xy
    \draw[->] (1.7,2.2) -> (-5.8,8.9);  % differential on t^-1xy
    \node at (-6,2) {$\cdots$};
    \node at (1.5,4) {$\cdots$};
    \node at (-7,3.5) {$\cdots$};
    \end{tikzpicture}
\end{figure}
%\end{comment}

To solve the extension problem we employ a strategy exactly analogous to the one used in \cite[Prop. 12]{SpeirsTruncated}.  The class $z_{(m,s)}$ in the homology of $B(m,s)$ gives a map of chain complexes 
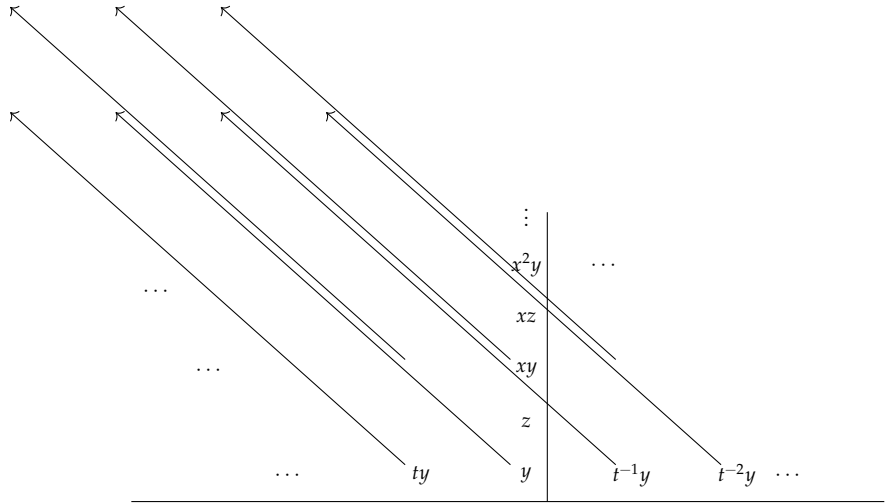
\[
f_{z_{(m,s)}} : \ZZ[m] \to \ZZ \otimes B(m,s)
\]
which may be promoted to a map of $\TT$-equivariant complexes, i.e.\ a map in the $\infty$-category $\D(\ZZ)^{B\TT}$. This category is equivalent to the $\infty$-category of mixed complexes, and so to check that $f_{z_{(m,s)}}$ is $\TT$-equivariant it suffices to check that the Connes operator on the taget $\ZZ \otimes B(m,s)$ is trivial. This follows from \cref{FinalClassificationLemma}. Now using the multiplicative $\TT$-equivariant map $\ZZ \to \THH(k)$ constructed in the proof of \cref{InfiniteCycleLemma}, we get a $\TT$-equivariant map
\[
\THH(k)[m] \to \THH(k) \otimes B(m,s) .
\]
The induced map on Tate spectral sequences now shows us that the $W(k)$-module $\pi_{2r}(\THH(k) \otimes B(m,s))^{t\TT}$ is cyclic and hence determined by its length (cf. \cite[Prop. 12]{SpeirsTruncated} for the details). The length is precisely what we have established from the $E^\infty$-page calculation.

 Truncating the spectral sequence we now see that
\[
\pi_{2r}(\THH(k) \otimes B(p^{v+1}m',p^us'))^{h\TT} \simeq \left\{
	\begin{array}{ll}
		  W_{v-u+2}(k) & \mbox{if } 2r \geq p^{v+1}m' \\
		  W_{v-u+1}(k) & \mbox{if } 2r < p^{v+1}m'
	\end{array}
\right.
\]
This completes the proof of \emph{(1)}.

The arguments in case \emph{(2)}, where  $m$ and $s$ are both odd, are very similar. % \todo{include these arguments?}
\end{proof}

\begin{remark}
In particular this shows that $\TP(A,B,I)$ is non-trivial. This may be contrasted with the Cuntz-Quillen result \cite{CuntzQuillenExcision}, that $\HP((A,B,I)/\QQ)$ is trivial, i.e.\ that rational periodic cyclic homology satisfies excision. Similarly it is shown in \cite{Horiuchi} that $\TP$ does not satisfy nil-invariance. Here again it is a result of Goodwillie \cite{Goodwillie1985} that rational periodic cyclic homology \emph{does} satisfy nil-invariance.
\end{remark}

Suppose $m'$ is even and let $t_{ev} = t_{ev}(p,r,m')$ be the unique positive integer such that $p^{t_{ev}-1}m' \leq 2r < p^{t_{ev}}m'$ (or zero if such $t_{ev}$ does not exist, i.e.\ if $m'$ is too big). Then we may restate the $\TC^-$ calculations as
\[
 \pi_{2r}((\THH(k) \otimes B(p^vm',p^us'))^{h\TT}) \simeq \left\{
        \begin{array}{ll}
           W_{v-u+1}(k)  & v < t_{ev}  \\
            W_{v-u}(k) & \quad  v \geq t_{ev}
        \end{array}
    \right. .
\]
Similarly if $m'$ and $s'$ are both odd, let $t_{od} = t_{od}(p,r,m')$ be the unique positive integer such that $p^{t_{od}-1}m' \leq 2r + 1 < p^{t_{od}}m'$ (or zero if such $t_{od}$ does not exist, i.e.\ if $m'$ is too big). Then
\[
 \pi_{2r+1}((\THH(k) \otimes B(p^vm',p^us'))^{h\TT}) \simeq \left\{
        \begin{array}{ll}
           W_{v-u+1}(k)  &  v < t_{od}  \\
            W_{v-u}(k) & \quad  v \geq t_{od}
        \end{array}
    \right. .
\]

\subsection{The case $p=2$}\label{Casep=2}

We now deal with the case when $k$ has characteristic two. From \cref{BiRelativeTHHp=2} we see that the case $m$ even and $s$ odd, is missing from \cref{TPandTCminus}. 

\begin{proposition}\label{TPp=2}
Assume $p=2$, $m$ even and $s$ odd. Write $m = 2^v m'$ with $m'$ odd, and let $s' \mid m'$. Then the homotopy groups of $(\THH(k) \otimes B(m,s))^{t\TT}$ and $(\THH(k) \otimes B(m,s))^{h\TT}$ are concentrated in odd degrees where they are isomorphic to $k$. 
\end{proposition}
\begin{proof}
When $v = 0$ we are in the case \emph{(2)} of \cref{TPandTCminus} so 
\[
\pi_*(\THH(k) \otimes B(m',s'))^{t\TT} = 0
\]
 and 
 \[
 \pi_*(\THH(k) \otimes B(m',s'))^{h\TT} = k \quad  \text{ when } \quad  * = 2r+1  \geq m'.
 \]
From the high co-connectivity of the Frobenius we conclude that 
\[
\pi_*(\THH(k) \otimes B(2m',s'))^{t\TT} = k
\]
 when $*$ is odd, and zero otherwise. Now the Tate spectral sequence 
\[
E^2 = k[t^{\pm 1}, x]\{z,w\} \quad \Rightarrow \quad \pi_*(\THH(k) \otimes B(2m',s'))^{t\TT}
\]
must collapse on the $E^2$-page, from which we conclude that $d^2(w) \overset{.}= (tx)z$. Truncating the spectral sequence we see the homotopy fixed point spectral sequence has $E^3 = E^\infty = (k[t,x]/(tx))\{z\}$ so 
\[
\pi_*(\THH(k) \otimes B(2m',s'))^{h\TT} = k
\]
in every odd degree. Now an induction argument shows that this pattern continues for all $v > 1$. 
% \todo{clean up this argument, and draw the Tate spectral sequence}
\end{proof}

%%%%%%%%%%%%%%%%%%%%%%%%%%%%%%%%%%%%%%%%
%%%%%%%%%%%%%%%%%%%%%%%%%%%%%%%%%%%%%%%%
%%%%%%%%%%%%%%%%%%%%%%%%%%%%%%%%%%%%%%%%

\section{Topological cyclic homology of coordinate axes}\label{TCsection}
\noindent
In this section we complete the proof of \cref{TCCoordinateAxes}, and thus of \cref{KtheoryCoordinateAxes}.
We start with the case where $p >2$.
For $1 \leq m'$ and $s' \mid m'$ let 
\[
\TP(m',s') := \prod_{0 \leq v} \prod_{u \leq v} \left( (\THH(k) \otimes B(p^vm',p^us'))^{t\TT} \right)^{\cyc_d(p^us')}
\]
and
\[
\TC^-(m',s') := \prod_{0 \leq v} \prod_{u \leq v} \left( (\THH(k) \otimes B(p^vm',p^us'))^{h\TT} \right)^{\cyc_d(p^us')}
\]
Thus, $\TC(A,B,I)$ is identified with the product of the equalizer of the maps $\phi_p, can : \TC^-(m',s') \rightarrow \TP(m',s')$ as $m'$ and $s'$ vary accordingly. Let us denote by $\TC(m',s')$ the equalizer
\[
\TC(m',s') \to \TC^-(m',s') \underset{can}{\overset{\phi_p}\rightrightarrows} \TP(m',s')
\]
A priori the homotopy groups of $\TC(m',s')$ sit in a long exact sequence with those of $\TC^-(m',s')$ and $\TP(m',s')$. However this sequence splits into short exact sequences since $\TC^-(m',s')$ and $\TP(m',s')$ are concentrated in either even or odd degrees, depending on the parity of $m'$ and $s'$, cf.\  \cref{TPandTCminus}.

%Suppose first that $p > 2$. 
If $m'$ and $s'$ are even then the Frobenius map
\[
\pi_{2r}(\THH(k) \otimes B(p^vm',p^us'))^{h\TT} \longto \pi_{2r}(\THH(k) \otimes B(p^{v+1}m',p^us'))^{t\TT}
\]
is an isomorphism for $0 \leq v < t_{ev}$, and the canonical map
\[
\pi_{2r}(\THH(k) \otimes B(p^vm',p^us'))^{h\TT} \longto \pi_{2r}(\THH(k) \otimes B(p^{v}m',p^us'))^{t\TT}
\]
is an isomorphism for $t_{ev} \leq v$. Thus we have a map of short exact sequences
%\[
%\begin{tikzcd}[row sep=0.8em, column sep=0.8em]
%\prod_{t_{ev} \leq v} \prod_{u \leq v} (W_{v-u}(k))^{\cyc_d(p^us')} \arrow[d, "\phi-can"] \arrow[r, %hook] & \TC^{-}_{2r}(m',s') \arrow[d, "\phi - can"] \arrow[r, two heads] & \prod_{0 \leq v < t_{ev}} %\prod_{u \leq v} (W_{v - u+1}(k))^{\cyc_d(p^us')} \arrow[d, "\overline{\phi - can}"]  \\
%\prod_{t_{ev} \leq v} \prod_{u \leq v} (W_{v-u}(k)))^{\cyc_d(p^us')} \arrow[r, hook] & \TP_{2r}(m',s') \arrow[r, two heads] & \prod_{0 \leq v < t_{ev}} \prod_{u \leq v} (W_{v-u}(k)))^{\cyc_d(p^us')}
%\end{tikzcd}
%\]
\[
\begin{tikzcd}[row sep=3em, column sep=1.4em]
\prod_{t_{ev} \leq v} \prod_{u \leq v} (W_{v-u}(k))^{\cyc_d(p^us')}  \arrow[r, "\phi_p-can"]   \arrow[d, hook]  &    
\prod_{t_{ev} \leq v} \prod_{u \leq v} (W_{v-u}(k)))^{\cyc_d(p^us')}  \arrow[d, hook]
\\
\TC^{-}_{2r}(m',s') \arrow[r, "\phi_p-can"]  \arrow[d, two heads]   &
 \TP_{2r}(m',s') \arrow[d, two heads]
\\
\prod_{0 \leq v < t_{ev}} \prod_{u \leq v} (W_{v - u+1}(k))^{\cyc_d(p^us')}  \arrow[r, "\overline{\phi_p-can}"]    &
\prod_{0 \leq v < t_{ev}} \prod_{u \leq v} (W_{v-u}(k)))^{\cyc_d(p^us')}
\end{tikzcd}
\]
The top horizontal map is an isomorphism, and the bottom horizontal map is an epimorphism so, by the snake lemma, we conclude that
\[
\TC_{2r}(m',s') \simeq \prod_{u \leq t_{ev}} (W_{t_{ev} - u}(k))^{\cyc_d(p^us')} .
\]

\noindent Now suppose $m'$ and $s'$ are odd. Then the Frobenius map
\[
\pi_{2r+1}(\THH(k) \otimes B(p^vm',p^us'))^{h\TT} \to \pi_{2r+1}(\THH(k) \otimes B(p^{v+1}m',p^us'))^{t\TT}
\]
is an isomorphism for $0 \leq v < t_{od}$, and the canonical map
\[
\pi_{2r+1}(\THH(k) \otimes B(p^vm',p^us'))^{h\TT} \to \pi_{2r+1}(\THH(k) \otimes B(p^{v}m',p^us'))^{t\TT}
\]
is an isomorphism for $t_{od} \leq v$. Thus we have a map of short exact sequences
%\[
%\begin{tikzcd}[row sep=1em, column sep=1em]
%\prod_{t_{od} \leq v} \prod_{u \leq v} (W_{v-u}(k))^{\cyc_d(p^us')} \arrow[d, "\phi-can"] \arrow[r, hook] & \TC^{-}_{2r+1}(m',s') \arrow[d, "\phi - can"] \arrow[r, two heads] & \prod_{0 \leq v < t_{od}} \prod_{u \leq v} (W_{v - u+1}(k))^{\cyc_d(p^us')} \arrow[d, "\overline{\phi - can}"]  \\
%\prod_{t_{od} \leq v} \prod_{u \leq v} (W_{v-u}(k)))^{\cyc_d(p^us')} \arrow[r, hook] & \TP_{2r+1}(m',s') \arrow[r, two heads] & \prod_{0 \leq v < t_{od}} \prod_{u \leq v} (W_{v-u}(k)))^{\cyc_d(p^us')}
%\end{tikzcd}
%\]
\[
\begin{tikzcd}[row sep=3em, column sep=1.4em]
\prod_{t_{od} \leq v} \prod_{u \leq v} (W_{v-u}(k))^{\cyc_d(p^us')}  \arrow[r, "\phi_p-can"]   \arrow[d, hook]  &    
\prod_{t_{od} \leq v} \prod_{u \leq v} (W_{v-u}(k)))^{\cyc_d(p^us')}  \arrow[d, hook]
\\
\TC^{-}_{2r+1}(m',s') \arrow[r, "\phi_p-can"]  \arrow[d, two heads]   &
 \TP_{2r+1}(m',s') \arrow[d, two heads]
\\
\prod_{0 \leq v < t_{od}} \prod_{u \leq v} (W_{v - u+1}(k))^{\cyc_d(p^us')}  \arrow[r, "\overline{\phi_p-can}"]    &
\prod_{0 \leq v < t_{od}} \prod_{u \leq v} (W_{v-u}(k)))^{\cyc_d(p^us')}
\end{tikzcd}
\]
The top horizontal map is an isomorphism, and the bottom horizontal map is en epimorphism so, by the snake lemma, we conclude that
\[
\TC_{2r+1}(m',s') \simeq \prod_{u \leq t_{od}} (W_{t_{od} - u}(k))^{\cyc_d(p^us')} .
\]
This finishes the proof of \cref{TCCoordinateAxes} in the case $p >2$.

\subsection{The case $p=2$}
If $p=2$ then we must correct slightly the definition of $\TP(m',s')$ and $\TC^-(m',s')$. Suppose $m'$ and $s'$ are odd. To deal with the case $m = p^vm'$  even and $s = p^us'$ even, we let
\[
\TP(m',s')_{ev} := \prod_{1 \leq v} \prod_{1 \leq u \leq v} \left( (\THH(k) \otimes B(2^vm',2^us'))^{t\TT} \right)^{\cyc_d(2^us')}
\]
and
\[
\TC^-(m',s')_{ev} := \prod_{1 \leq v} \prod_{1 \leq u \leq v} \left( (\THH(k) \otimes B(2^vm',2^us'))^{h\TT} \right)^{\cyc_d(2^us')}
\]
Both of these spectra are concentrated in even degrees, with their homotopy groups given by \cref{TPandTCminus}. As a result we see that, for $m' \geq 1$ odd, and $s' \mid m'$ odd, 
\[
\TC_{2r}(m',s')_{ev} \simeq \prod_{1 \leq u \leq t_{ev}} W_{t_{ev}-u}(k)^{\oplus \cyc_d(2^us')} .
\]

To deal with the case where $m$ is even and $s \mid m$ is odd, let
\[
\TP(m',s')_{ev,od} := \prod_{0 \leq v} \left( (\THH(k) \otimes B(2^vm',s'))^{t\TT} \right)^{\cyc_d(s')}
\]
and
\[
\TC^-(m',s')_{ev,od} := \prod_{0 \leq v} \left( (\THH(k) \otimes B(2^vm',s'))^{h\TT} \right)^{\cyc_d(s')}
\]
(note that when $v = 0$, we have $m = m'$ odd, but we must include this case since the Frobenius connects it with the case $v=1$)
By \cref{TPp=2}, $\TP(m',s')_{ev,od}$ and $\TC^-(m',s')_{ev,od}$ are concentrated in odd degrees, where they are isomorphic to $k$. Thus
\[
\TC_{2r+1}(m',s')_{ev,od} \simeq \prod_{0 \leq v \leq t_{ev}} k^{\oplus \cyc_d(s')} .
\]
This completes the proof of \cref{TCCoordinateAxes} when $p=2$.

%%%%%%%%%%%%%%%%%%%%%%%%%%%%%%%%%%%%%%%%
%%%%%%%%%%%%%%%%%%%%%%%%%%%%%%%%%%%%%%%%
%%%%%%%%%%%%%%%%%%%%%%%%%%%%%%%%%%%%%%%%

\section{The characteristic zero case}\label{CharZeroSection} 
\noindent
In this section we compute the relative cyclic homology and the bi-relative $K$-theory of $A_d = k[x_1,\dots , x_d]/(x_ix_j)_{i \ne j}$ in the case that $k$ is an ind-smooth $\QQ$-algebra. In the case where $k/\QQ$ is a field extension, the results in this section were found, by different means, already in 1989 by Geller, Reid and Weibel \cite[Theorem 7.1.]{GellerReidWeibel}. %In the following we write $\HH(A) = \HH(A/\ZZ)$, i.e.\ with integral base, unless otherwise indicated.

We proceed as in \cite[Section 3.9]{hesselholt2005k}. We will compute the relative groups 
\[
\HC_q((A, I)/\ZZ) \otimes \QQ \simeq \HC_q((A,I) /\QQ)
\] using our understanding of the $\TT$-homotopy type of the cyclic bar construction for $\Pi^d$, as found in \cref{FinalClassificationLemma}.
First we need the following analogue of \cref{SplittingLemma}

\begin{lemma}
Let $k$ be a ring, $\Pi$ a pointed monoid, and $k(\Pi)$ the pointed monoid algebra. There is a natural equivalence of $\TT$-spectra
\[
\HH(k(\Pi)) \overset{\sim}\leftarrow \HH(k) \otimes \Bcy(\Pi)
\]
\end{lemma}
\begin{comment}
\begin{proof}
\todo[line]{write a proof. This may also be needed in the paper on truncated polynomial algebras?}
\end{proof}
\end{comment}

So the arguments from \cref{IndexLemma} still work, yielding a description of the relative (and bi-relative) Hochshild homology spectrum of $(A_d,I_d)$ (and $(A_d,B_d,I_d)$). First we note that the relative Hochschild homology differs only slightly from the absolute version. The sole difference is that we ``cut out'' the part of the space $\Bcy(\Pi^d)$ given by $\Bcy(\Pi^d; \overline{\emptyset})$, i.e.\ the part corresponding to the unique cyclic word of length zero. %Expressed in the form, as before we have the following cofiber sequence
Since the spaces $B(m,s)$ for $m$ even and $s$ odd have torsion integral homology they disappear in the rational case. So we conclude that
\begin{eqnarray*}
 \HH(A,I) \simeq &&\bigoplus_{\substack{m \geq 2 \\ even} } \bigoplus_{\substack{s \mid m \\ even}} \left( \HH(k) \otimes B(m,s)  \right)^{\oplus \cyc_d(s)} 
\\
 \oplus &&
  \bigoplus_{\substack{m \geq 2 \\ odd}} \bigoplus_{\substack{ s \mid m \\  odd}} \left( \HH(k) \otimes B(m,s)  \right)^{\oplus \cyc_d(s)} 
  \\
\oplus &&   \bigoplus_{i\geq 1} \left( \HH(k) \otimes (\TT/C_{i})_+ \right)^{\oplus d}
\end{eqnarray*}
where $B(m,s)$ is given by \cref{FinalClassificationLemma}.
The bottom summands, with terms $\HH(k) \otimes (\TT/C_i)_+$, corresponds to the spaces $\Bcy(\Pi^d, \overline{x_1^i})$ , $\Bcy(\Pi^d, \overline{x_2^i})$ etc. These bottom summands disappear when looking at the bi-relative theory, $\HH(A,B,I)$.

\begin{theorem} \label{CyclicHomology}
Let $k$ be any commutative unital $\QQ$-algebra. There is an isomorphism
\[
\HC_0(A_d/\QQ) \simeq A_d
\]
and, for $q \geq 1$, an isomorphism
\begin{eqnarray*}
\HC_q((A_d, I_d)/\QQ) \simeq && \bigoplus_{\substack{ m \geq 2 \\ even}} \bigoplus_{\substack{s \mid m \\ even}} \left( \HH_{q+1-m}(k) \ \otimes \QQ \right)^{\oplus \cyc_d(s)} 
\\
\oplus &&  \bigoplus_{\substack{m \geq 2 \\ odd}} \bigoplus_{\substack{s \mid m \\ odd}} \left( \HH_{q +1 - m}(k) \otimes \QQ \right)^{\oplus \cyc_d(s)} 
\\
\oplus && \bigoplus_{i \geq 1} \left( \HH_q(k) \otimes \QQ \right)^{\oplus d}
\end{eqnarray*}
\end{theorem}
\begin{proof}
We must compute $\pi_q(X_{h\TT})$ where $X$ ranges over the summands in the above decomposition of $\Bcy(\Pi^d)$. If $X = \HH(k) \otimes S^{\lambda_j} \otimes (\TT/C_i)_+$ for some complex $\TT$-representation $\lambda_j$ of complex dimension $j$, then since,
\[
\left( \HH(k) \otimes S^{\lambda_j} \otimes (\TT/C_i)_+ \right)_{h\TT} \simeq \left( \HH(k) \otimes S^{\lambda_j} \right)_{hC_i} ,
\] 
we may regard the $C_i$-homology spectral sequence
\[
E^2_{s,t} = H_s(C_i , \pi_t(\HH(k) \otimes S^{\lambda_j}) \otimes \QQ) \Rightarrow \pi_{s+t}(\HH(k) \otimes S^{\lambda_j})_{hC_i} \otimes \QQ
\]
Since the rational group homology of $C_i$ is concentrated in degree $0$, the edge homomorphism 
\[
H_0(C_i, \pi_q(\HH(k) \otimes S^{\lambda_j}) \otimes \QQ) \overset{\simeq}{\longto} \pi_q\left( \left( \HH(k) \otimes S^{\lambda_j} \right)_{hC_i} \right) \otimes \QQ
\]
is an isomorphism. Furthermore, since the $C_i$-action on  $\HH(k) \otimes S^{\lambda_j}$ extends to a $\TT$-action, the induced action on homotopy groups is trivial. We conclude that
\[
\pi_q\left( \left( \HH(k) \otimes S^{\lambda_j} \otimes (\TT/C_i)_+ \right)_{h\TT} \right) \cong \HH_{q - 2j}(k)
\]
The result follows.
\end{proof}

In the case where $k$ is an algebraic field extension of $\QQ$ \cref{CyclicHomology} was established by Geller, Reid and Weibel \cite[Example 5.7.]{GellerReidWeibel}. Indeed, then $k$ is ind-smooth over $\QQ$ (cf.\ \cite[07BV]{StacksProject})  so the Hochschild-Kostant-Rosenberg theorem calculates the Hochschild homology in terms of differential forms. Similarly, the following theorem is proved in  \cite[Theorem 7.1.]{GellerReidWeibel} when $k$ is a field extension of $\QQ$. We use their counting function $c_{d-1}(q)$ which we recall in \cref{CyclicAppendix}, \cref{GRWfunction}.

\begin{theorem}
Suppose $k$ is an ind-smooth $\QQ$-algebra. Let $d \geq 2$ and consider the ring $A_d = k[x_1, \dots , x_d]/(x_ix_j)_{i \ne j}$, and let $I_d = (x_1, \dots , x_d)$. Then
\[
K_q(A_d,I_d)  \cong   k^{\oplus c_{d-1}(q)} \oplus (\Omega^1_{k/\QQ})^{\oplus c_{d-1}(q-1)} \oplus \dots \oplus (\Omega^{q-2}_{k/\QQ})^{\oplus c_{d-1}(2)} .
\]
\end{theorem}
\begin{proof}
By \cite[Corollary 0.2]{CortinasKABI} we have $K_n(A_d,I_d) \simeq \HC_{n-1}(A_d,B_d,I_d)$. Here we use that $K_q(A_d,I_d)$ is a rational vector space \cite[Theorem 0.1]{CortinasKABI}, so no further rationalization is necessary. By \cref{CyclicHomology} we are reduced to Hochschild homology calculations. By the Hochschild-Kostant-Rosenberg theorem \cite[Theorem 9.4.7.]{WeibelBook} we have
\[
\HH_n(k/\QQ) \simeq \Omega^{n}_{k/\QQ}.
\]
Thus, 
\begin{eqnarray*}
K_q(A_d,I_d) &\simeq & \HC_{q-1}(A_d,B_d,I_d)
\\ &\simeq & k^{\oplus c_{d-1}(q)} \oplus (\Omega^1_{k/\QQ})^{\oplus c_{d-1}(q-1)} \oplus \dots \oplus (\Omega^{q-2}_{k/\QQ})^{\oplus c_{d-1}(2)}
\end{eqnarray*}
This completes the proof.
\end{proof}

\begin{remark}
If $k$ is a field extension of $\QQ$ for which we know the transcendence degree of $k$ over $\QQ$ then this result completely determines the relative $K$-theory, since $\dim_k \Omega^1_{k/\QQ} = tr.deg(k/\QQ)$. For example if $k/\QQ$ is algebraic then $\Omega^1_{k/\QQ} = 0$ and so the calculation reduces to
\[
K_q(A_d,I_d) \simeq k^{\oplus c_{d-1}(q)}.
\]
%Sketch of proof that $\dim_k \Omega^1_{k/\QQ} = tr.deg(k/\QQ)$: Suppose first that $tr.deg(k/\QQ) = 1$ i.e. $k = \QQ(t)$. Then $\Omega^1_{\QQ(t)/\QQ}$ is the localization of $\Omega_{\QQ[t]/\QQ} \simeq \QQ[t]dt$ which is just $\QQ(t)dt$. For general $k$ let $S$ be a transcendence basis and $K = \QQ(S)$ then have $\QQ \subs K \subs k$ with $k/K$ algebraic. So $\Omega^1_{k/K} = 0$ and by the Jacobi-Zariske sequence $k \otimes_K \Omega^1_{K\QQ}$ injects into $\Omega^1_{k/\QQ}$. Now one reduced to the case where $k/K$ is a finite extension. Here one uses trickery with the conormal sequence. The general case follows since $\Omega^1$ commutes with colimits. 
\end{remark}

%%%%%%%%%%%%%%%%%%%%%%%%%%%%%%%%%%%%%%%%
%%%%%%%%%%%%%%%%%%%%%%%%%%%%%%%%%%%%%%%%
%%%%%%%%%%%%%%%%%%%%%%%%%%%%%%%%%%%%%%%%

\section{Appendix: counting cyclic words}\label{CyclicAppendix}
\noindent
In this section we use some counting techniques inspired by \cite{GellerReidWeibel}. We first deal with all words, then with cyclic words. 
Let $A_d(m)$ denote the number of words in $d$ letters of length $m$ having no cyclic repetitions.

\begin{lemma}\label{CountingLemma}
For all $d \geq 1$ and all $m \geq 1$ we have 
\[
A_d(m) = (d-1)^m + (-1)^m (d-1).
\]
\end{lemma}
\begin{proof}
Define the auxiliary counter, $B_d(m)$ counting words $\omega = w_1w_2 \dots w_m$ with no allowed repetitions, except that we require $w_m = w_1$.  Then
\[
A_d(m) + B_d(m) = d(d-1)^{m-1}.
\]
To see this, consider a word $\omega = w_1w_2 \dots w_m$. There are $d$ choices for $w_1$, and $d-1$ choices for $w_2, w_3, \dots$ and $ w_{m-1}$. Finally for $w_m$ there are again $d-1$ choices since choosing any letter different from $w_{m-1}$ gives either a type $A$ ($d-2$ possibly choices) or a type $B$ word (one choice, namely $w_m = w_1$). 
The formula for $A_d(m)$ now follows by induction using the equation $B_d(m) = A_d(m-1)$ (for $m \geq 2$). This last equality is true since deleting the last letter of a type $B$ word yields a word with no cyclic repetitions.
\end{proof}

Let $\widetilde{\cyc}_d(s)$ denote the number of words in $d$ letters of length $s$, period $s$ and having no cyclic repetitions (see \cref{NoReps}) and let $\cyc_d(s)$ denote the number of \emph{cyclic} words in $d$ letters of length $s$, period $s$ and having no cyclic repetitions. Then $\cyc_d(s) = \frac{1}{s}\widetilde{\cyc}_d(s)$. Also we have
\[
A_d(m) = \sum_{s \mid m}\widetilde{\cyc}_d(s)
\]
So using M\"obius inversion we have
\[
\widetilde{\cyc}_d(s) = \sum_{j \mid s} \mu(\frac{s}{j})A_d(j).
\]
hence we obtain a formula for $\cyc_d(s)$,
\[
\cyc_d(s) = \frac{1}{s} \sum_{j \mid s} \mu(\frac{s}{j}) ((d-1)^j + (-1)^j(d-1)).
\]
Below is a table of the first few values.

\begin{center}
\setlength{\extrarowheight}{5pt}
 \begin{tabular}{|c|  c c c|}
 \hline
  s & $\cyc_d(s)$  & $d=3$  & $d=4$ \\ [0.7ex] 
 \hline
 1 & 0 & 0  & 0  \\ 
 \hline
 2 & $\frac{1}{2}((d-1)^2 + (d-1))$  & 3 & 6  \\ [1ex] 
 \hline
 3 & $\frac{1}{3}((d-1)^3 - (d-1))$  &  2 & 8 \\ [1ex] 
 \hline
 4 & $\frac{1}{4}((d-1)^4 - (d-1)^2)$ & 3 & 18  \\ [1ex] 
 \hline
 5 & $\frac{1}{5}((d-1)^5 - (d-1))$ & 6  & 48 \\ [1ex] 
 \hline
  6 & $\frac{1}{6}((d-1)^6 - (d-1)^3 - (d-1)^2 + (d-1))$ & 9 & 116  \\ [1ex] 
 \hline
   7 & $\frac{1}{7}((d-1)^7 - (d-1))$ & 18 & 312  \\ [1ex] 
 \hline
   8 & $\frac{1}{8}((d-1)^8 - (d-1)^4)$ & 30  & 810 \\ [1ex] 
 \hline
   9 & $\frac{1}{9}((d-1)^9 - (d-1)^3)$ & 56 & 2184  \\ [1ex] 
    \hline
   10 &$\frac{1}{10}((d-1)^{10} - (d-1)^5 - (d-1)^2 + (d-1))$   & 99 &  5880 \\ [1ex] 
       \hline
   11 & $\frac{1}{11}((d-1)^{11} - (d-1))$ & 186 &  16104 \\ [1ex] 
       \hline
   12 & $\frac{1}{12}((d-1)^{12} - (d-1)^6 - (d-1)^4 + (d-1)^2)$ & 335 & 44220  \\ [1ex] 
 \hline
\end{tabular}
\end{center}
If we want to count the number of cyclic words without cyclic repetitions, with length $m$ but with no restrictions on period then we can just sum $\cyc_d(s)$ over all divisors $s$ of $m$. For example there are $3 + 2 + 9 = 14$ cyclic words in $d=3$ letter, with no repetitions having length $m = 6$.

We may describe the function $c_{d-1}(m)$ that Geller, Reid, and Weibel introduce as follows
\begin{equation}\label{GRWfunction}
c_{d-1}(m) = \sum_{\substack{s \mid m \\ s {\equiv} m\mod{2}}} \cyc_{d}(s) .
\end{equation}
That is, $c_{d-1}(m)$ is the number of cyclic words without repetitions, having length $m$ and period $s$ where $s$ and $m$ have the same parity.

%%%%%%%%%%%%%%%%%%%%%%%%%%%%%%%%%%%%%%%%
%%%%%%%%%%%%%%%%%%%%%%%%%%%%%%%%%%%%%%%%
%%%%%%%%%%%%%%%%%%%%%%%%%%%%%%%%%%%%%%%%
\begin{comment}
\section{Appendix: Examples}\label{Examples}
\noindent

The following tables provide an overview of the specific groups$ K_q(A_d,I_d)$ for small choices of $p,q$ and $d$.

\begin{center}
For $p=2$
\setlength{\extrarowheight}{5pt}
 \begin{tabular}{|c| c c c|}
  \hline
  $q$ & $d=2$  & $d=3$  & $d=4$ \\ [0.7ex] 
 \hline
 1 & 0 & 0  & 0  \\ 
 \hline
 2 & $\frac{1}{2}((d-1)^2 + (d-1))$  & 3 & 6  \\ [1ex] 
 \hline
 3 & $\frac{1}{3}((d-1)^3 - (d-1))$  &  2 & 8 \\ [1ex] 
 \hline
 4 & $\frac{1}{4}((d-1)^4 - (d-1)^2)$ & 3 & 18  \\ [1ex] 
 \hline
 5 & $\frac{1}{5}((d-1)^5 - (d-1))$ & 6  & 48 \\ [1ex] 
  \hline
\end{tabular}
\end{center}

\end{comment}

%\nocite{*}
\bibliographystyle{siam}
\bibliography{bibfileNewA}

\end{document}